\documentstyle[amstex,11pt,theorem,amscd]{amsart}

\edef\oldatcatcode{\the\catcode`\@}\catcode`\@=11
\@mparswitchfalse
 \newcommand{\NOTE}[1]{\ifinner\else
   {\marginpar{\hfill\smash{\tiny\normalshape\mediumseries #1}}}\fi}
\renewcommand{\rom}[1]{{\normalshape#1}}
\let\OLDref=\ref
\newcommand{\thmrom}{}
\renewcommand{\ref}[1]{\thmrom{\OLDref{#1}}}
\def\eqref#1{\rom{\tagform@{\ref{#1}}}}
\newenvironment{Proof}{\begin{ProofwCaption}{Proof}}{\end{ProofwCaption}}
\newenvironment{ProofwCaption}[1]%
  {\addvspace\theorempreskipamount \noindent{\it #1.}\rm}%
  {\qed \par \addvspace\theorempostskipamount}
\newcommand{\noqed}{\unskip\def\qed{}} 
\catcode`\@=\oldatcatcode

\def\PP{{\Bbb P}}
\def\ZZ{{\Bbb Z}}

\newcounter{thm}
\theoremstyle{plain}
\newtheorem{Theorem}[thm]{Theorem}
\newtheorem{Lemma}[thm]{Lemma}
\newtheorem{Proposition}[thm]{Proposition}
\newtheorem{Corollary}[thm]{Corollary}
\newtheorem{Claim}{Claim}
\theoremstyle{definition}
\newtheorem{Definition}[thm]{Definition}
\theoremstyle{remark}
\newtheorem{Remark}[thm]{Remark}
\newtheorem{Remarks}[thm]{Remarks}
\newcommand{\dashto}{\mathrel{\dashrightarrow}}
\newsymbol\boxtimes      1202
\newsymbol\rtimes        226F
\newsymbol\sm 2272  
\newsymbol\shortmid 2370
\newcommand{\tr}[1]{\vphantom{#1}^{{\normalshape t}}\mskip-1.5mu #1}
\newcommand{\Chi}{{\cal X}}
\newcommand{\iso}{\cong}
\begin{document}

\title{The geometry of bielliptic surfaces in $\PP^4$}
\author[A. Aure, W. Decker, K. Hulek, S. Popescu, K. Ranestad]
{Alf Aure \and Wolfram Decker \and  Klaus Hulek, \\
Sorin Popescu \and  Kristian Ranestad \\ }

\maketitle

In 1988 Serrano \cite{Ser}, using Reider's method, discovered a minimal
bielliptic surface in $\PP^4$. Actually he showed that there is a unique
family of such surfaces and that they have degree 10 and sectional genus 6.
It is easy to see that the only other smooth surfaces with these invariants
are minimal abelian. There is a unique family of minimal abelian surfaces
in $\PP^4$; these arise as the smooth zero-schemes of sections of the
Horrocks-Mumford bundle \cite{HM}. The geometry of the Horrocks-Mumford
surfaces (smooth or not) has been intensively studied (see e.g. \cite{BHM}).

In this paper we shall, among other things, describe the geometry of the
embedding of the minimal bielliptic surfaces. A consequence of this description
will be the existence of smooth nonminimal bielliptic surfaces of degree
15 in $\PP^4$.

The starting point of our investigations is the following consequence
of Serrano's result: The bielliptic surface of degree 10 has a fibration
onto an elliptic curve whose fibres are all plane cubic curves. If $V$
is the union of the planes of these plane cubics, then $V$ has degree 5
and is the union of the trisecants to an elliptic quintic scroll in
$\PP^4$ with the scroll as its double surface. To see this assume that
two members of the family of plane cubic curves span only a hyperplane,
then the residual curve in that hyperplane section is a rational quartic
curve, but that is impossible since it would dominate the elliptic base
of the pencil. Thus in the dual space, the planes of $V$ correspond to
the lines of a smooth elliptic scroll. A first account of $V$ was
given by Segre \cite{Seg}.

In a recent paper Catanese and Ciliberto \cite{CC} computed the cohomology
of antipluricanonical divisors on the symmetric products of an elliptic
curve using Heisenberg invariants. This lead us to see that the
minimal desingularization $\widetilde V$ of $V$ contains 8 bielliptic surfaces.
These are embedded in $\PP^4$ via the map $\widetilde V \to V$. In order to
understand this map well it is useful to interpret $V$ as the set of all
singular quadrics through a quintic elliptic normal curve $E$ in $\PP^4$.
Blowing up the $\PP^2$-bundle $\widetilde V$ in a certain section we get a
$\PP^1$-bundle $\widetilde W$ over $S^2 E$. It turns out that $\widetilde W$ is
the
natural desingularization of the secant variety $W$ of the quintic curve
$E$. The threefold $\widetilde W$ contains 8 bielliptic surfaces blown up in 25
points. Via the map $\widetilde W \to W$ these surfaces are embedded as smooth
surfaces of degree 15. In this way we can also find nonminimal abelian
surfaces of degree 15 in $\PP^4$ whose minimal model is isogeneous to a
product. These surfaces are $(5,5)$-linked to Horrocks-Mumford surfaces,
their existence being known before. On the other hand the nonminimal bielliptic
surfaces lie on a unique quintic, namely $W$.

The above duality between $V$ and the elliptic scroll can also be interpreted
in terms of the quadro-cubic Cremona transformation of Semple \cite{Sem},
\cite{SR}. We explain how to construct the degree 15 surfaces with the
help of this Cremona transformation.

In a subsequent paper we will study the syzygies and thus also the
degenerations of abelian and bielliptic surfaces in $\PP^4$. A consequence
will be the existence of a family of smooth nonminimal abelian surfaces of
degree 15 lying on only one quintic hypersurface.

We finally remark that the quintic elliptic scroll and the abelian and
bielliptic surfaces of degree 10 and 15 are essentially the only smooth
irregular surfaces known in $\PP^4$ (all others can be derived via finite
morphisms  $\PP^4\to\PP^4 $).

The authors thank the DFG for financial support in the framework of the
{\em Schwer\-punkt\-pro\-gramm
``Kom\-plexe Man\-nig\-fal\-tig\-kei\-ten''}.

\setcounter{section}{-1}
\section{Heisenberg invariants on $\Bbb P^2$}

Here we collect some well-known facts about invariants of the
Schr\"odinger representation of $H_3$, the Heisenberg group of level
$3$. Let $x_0,x_1,x_2$ be a basis of
$\mathrm H^\circ(\cal O_{\Bbb P^2}(1))$ and consider the dual of the
Schr\"odinger representation of $H_3$ on
$V=\mathrm H^\circ(\cal O_{\Bbb P^2}(1))$ given by
\begin{equation}\label{(1)}
  \begin{aligned}
    \sigma_3(x_i)&=x_{i-1}\\
    \tau_3(x_i)&=\varepsilon_3^{-i}x_i\quad(\varepsilon_3=e^{2\pi i/3})
  \end{aligned}
\end{equation}
where $i$ is counted modulo 3 and $\sigma_3$ and $\tau_3$ generate
$H_3$. Note that
\begin{equation}\label{(1a)}
  [\sigma_3,\tau_3]=\varepsilon_3^{-1}\cdot\operatorname{id},
\end{equation}
hence $H_3$ is a central extension
\[
    1
    \to
    \mu_3
    \to
    H_3
    \to
    {\Bbb Z}_3\times{\Bbb Z}_3
    \to
    1.
  \]
The induced representation on
$\mathrm H^\circ(\cal O_{\Bbb P_2}(3))$ decomposes
into characters since $\sigma_3$ and  $\tau_3$ commute on the third
symmetric power of
$\mathrm H^\circ(\cal O_{\Bbb P^2}(1))$. By $(a,b)$ we denote
the character where $\sigma_3$ (resp.\ $\tau_3$) acts by
$\varepsilon_3^a$ (resp.\ $\varepsilon_3^b$)). Here again $a,b$ have to
be taken modulo 3. There is a pencil of invariant polynomials, called
the Hesse pencil, spanned by
\[
  x_0^3+x_1^3+x_2^3,\quad x_0x_1x_2
\]
and eight invariant polynomials corresponding to the eight non-trivial
characters:
\begin{align*}
  F_{(1,0)}&\colon\quad x_0^3+\varepsilon_3x_1^3+\varepsilon_3^2x_2^3\\
  F_{(2,0)}&\colon\quad x_0^3+\varepsilon_3^2x_1^3+\varepsilon_3x_2^3\\
  F_{(0,1)}&\colon\quad x_0x_1^2+x_1x_2^2+x_2x_0^2\\
  F_{(1,1)}&\colon\quad
x_0x_1^2+\varepsilon_3x_1x_2^2+\varepsilon_3^2x_2x_0^2\\
  F_{(2,1)}&\colon\quad
x_0x_1^2+\varepsilon_3^2x_1x_2^2+\varepsilon_3x_2x_0^2\\
  F_{(0,2)}&\colon\quad x_0^2x_1+x_1^2x_2+x_2^2x_0\\
  F_{(1,2)}&\colon\quad
x_0^2x_1+\varepsilon_3x_1^2x_2+\varepsilon_3^2x_2^2x_0\\
  F_{(2,2)}&\colon\quad x_0x_1^2+\varepsilon_3^2x_1x_2^2+\varepsilon_3x_2x_0^2
\end{align*}
where $F_{(a,b)}$ denotes the curve defined by the corresponding polynomial.
On each smooth member of the Hesse pencil the group acts by translation
of $3$-torsion points. There are precisely four singular members, namely:
\begin{align*}
  T_{(0,1)}&\colon\quad x_0x_1x_2\\
  T_{(1,1)}&\colon\quad(x_0+\varepsilon_3^2x_1+\varepsilon_3^2x_2)
                       (x_0+x_1+\varepsilon_3x_2)
                       (x_0+\varepsilon_3x_1+x_2)\\
  T_{(1,0)}&\colon\quad(x_0+\varepsilon_3x_1+\varepsilon_3^2x_2)
                       (x_0+\varepsilon_3^2x_1+\varepsilon_3x_2)
                       (x_0+x_1+x_2)\\
  T_{(1,2)}&\colon\quad(x_0+x_1+\varepsilon_3^2x_2)
                       (x_0+\varepsilon_3x_1+\varepsilon_3x_2)
                       (x_0+\varepsilon_3^2x_1+x_2)
\end{align*}
which equal $x_0^3+x_1^3+x_2^3+\lambda x_0x_1x_3$ for
$\lambda=\infty,-3\varepsilon_3^2,-3,-3\varepsilon_3$. For each $(i,j)$ the
subgroup of order 3 which is generated by $\sigma_3^i\tau_3^j$ fixes
the vertices of the triangle $T{(i,j)}$.

Consider the involution
\[
  \iota_3\colon(x_0,x_1,x_2)\mapsto(x_0,x_2,x_1).
\]
This involution leaves each member of the Hesse pencil invariant.
In fact choosing $(0,1,-1)$ as the origin it acts as $x\mapsto-x$ on
smooth members. The nontrivial characters come in pairs since
$\iota_3(F_{(a,b)})=F_{(-a,-b)}$.

\begin{Lemma}\label{L1}
  The curves $F_{(a,b)}$ are Fermat curves, and $H_3$ acts on each of them
  with translation by a 3-torsion point and multiplication by $\varepsilon_3$.
\end{Lemma}

\begin{Proof}
Let $\eta^3=\varepsilon_3$, $\mu^3=\frac19$. Then
\begin{align*}
  x_0^3+\varepsilon_3x_1^3+\varepsilon_3^2x_2^3
    & = x_0^3+(\eta x_1)^3+(\eta^2x_2)^3\\
  x_0^2x_1+x_1^2x_2+x_2^2x_0
    & = (\eta\mu(x_0+\varepsilon_3^2x_1+\varepsilon_3x_2))^3
        +(\eta^2\mu(x_0+\varepsilon_3x_1+\varepsilon_3^2x_2))^3\\
 &\qquad+(\mu(x_0+x_1+x_2))^3\\
  x_0^2x_1+\varepsilon_3x_1^2x_3+\varepsilon_3^2x_2^2x_0
    & = (\eta\mu(x_0+\varepsilon_3^2x_1+\varepsilon_3^2x_2))^3
        +(\eta^2\mu(x_0+\varepsilon_3x_1+x_2))^3\\
 &\qquad+(\mu(x_0+x_1+\varepsilon_3x_2))^3\\
  x_0^2x_1+\varepsilon_3^2x_1^2x_2+\varepsilon_3x_2^2x_0
    & = (\eta\mu(x_0+\varepsilon_3^2x_1+x_2))^3
        +(\eta^2\mu(x_0+\varepsilon_3x_1+\varepsilon_3x_2))^3\\
 &\qquad+(\mu(x_0+x_1+\varepsilon_3^2x_2))^3.
\end{align*}
Since $\iota_3(F_{(a,b)})=F_{(-a,-b)}$ all curves are Fermat curves.

The Fermat curve $F_{(a,b)}$ intersects each of the triangles $T_{(i,j)}$,
$(i,j)\ne\pm(a,b)$ in its vertices so $H_3$ has three subgroups of
order 3 with 3 fixed points on $F_{(a,b)}$. The fourth subgroup has
no fixed points, hence $H_3$ acts as stated.
\end{Proof}

\section{Threefolds containing bielliptic surfaces}\label{Par1}

In this part we will construct a $\Bbb P^2$-bundle over an elliptic curve
$E$, and a $\Bbb P^1$-bundle over the symmetric product $S^2E$ of the
elliptic curve containing bielliptic surfaces.

Choose a smooth element of the Hesse pencil
\[
  E=E_\lambda=\{x_0^3+x_1^3+x_2^3+\lambda x_0x_1x_2=0\}
\]
where $\lambda\ne\infty,-3,-3\varepsilon_3,-3\varepsilon_3^2$. We choose the
inflection point $p_0=(0,1,-1)$ to be the origin of $E$.

Let $\xi_0,\xi_1,\xi_2$ be a dual basis of
$x_0,x_1,x_2\in V=\Gamma(\cal O_{\Bbb P^2}(1))$. The induced action
of $H_3$ is given by
\begin{equation}\label{(2)}
  \begin{aligned}
    \sigma_3(\xi_i) &= \xi_{i-1}\\
    \tau_3(\xi_i)   &= \varepsilon_3^i\xi_i.
  \end{aligned}
\end{equation}

Note that in this case \eqref{(2)} implies
\begin{equation}\label{(2a)}
  [\sigma_3,\tau_3]=\varepsilon_3\cdot\operatorname{id}.
\end{equation}

Next consider the line bundle $\cal O_E(15p_0)$. Let $y_0,\ldots,y_{14}$ be
a basis of $\mathrm H^\circ(\cal O_E(15p_0))$ such that $H_{15}$, the
Heisenberg group of level 15, acts in the standard way, i.e., by
\begin{equation}\label{(3)}
  \begin{aligned}
    \sigma_{15}(y_i) &= y_{i-1}\\
    \tau_{15}(y_i)   &= \varepsilon_{15}^{-i}y_i
      \quad (\varepsilon_{15}=e^{2\pi i/15}).
  \end{aligned}
\end{equation}
{}From \eqref{(3)} it follows that
\begin{equation}\label{(3a)}
  [\sigma_{15}^5,\tau_{15}^5]=\varepsilon_{15}^{-10}\cdot\operatorname{id}
  =\varepsilon_3\cdot\operatorname{id}.
\end{equation}
Hence identifying $\sigma_{15}^5$ with $\sigma_3$ and $\tau_{15}^5$
with $\tau_3$ we get an isomorphism of the subgroup of $H_{15}$
generated by $\sigma_{15}^5$ and $\tau_{15}^5$ with
$H_3\subset\operatorname{SL}(V\spcheck)$, where the latter inclusion
is given by the Schr\"odinger representation. Since $y_0,\ldots,y_{14}$
generate $\cal O_E(15p_0)$ this gives an action of $H_3$ on the line bundle
$\cal O_E(15p_0)$ itself.

Hence we can consider the natural action of $H_3\times H_3$ on the rank 3
bundle $W_E=\cal O_E(15p_0)\otimes V$. Let $\Delta$ be the diagonal
of $H_3\times H_3$. Then $\Delta\iso H_3$ and
\begin{align*}
  \sigma_3(y_i\otimes x_j) &= y_{i-5}\otimes x_{j-1}\\
  \tau_3(y_i\otimes x_j)   &= \varepsilon_3^{-i-j}y_i\otimes x_j.
\end{align*}
It follows from \eqref{(1a)} and \eqref{(3a)} that the centre of $H_3$
acts trivially on $W_E$. Hence the quotient
\[
  \cal E_E=W_E\spcheck/\Delta
\]
is a rank 3 vector bundle over
\[
  E/\Bbb Z_3\times\Bbb Z_3=E.
\]

\begin{Lemma}\label{L2}
  \rom{(i)}$\cal E_E$ is stable of degree -5.
  \par\noindent\rom{(ii)}$\det\cal E_E=\cal O_E(-5p_0)$
\end{Lemma}

\begin{Proof}
(i) Since $\deg W_E=45$, the degree of $\cal E_E$ is clearly -5.
Now assume that $\cal F\subset\cal E_E\spcheck$ is a subbundle of rank
$r$ ($r=1,2$) and degree $d$ contradicting semistability, i.e., $d/r>5/3$.
Then ${\cal F}$ pulls back to a subbundle $\cal F'\subset W_E$ of degree
$9d>15r$.
This implies that $\cal F'\otimes\cal O_E(-15p_0)$ and hence
$W_E\otimes\cal O_E(-15p_0)$ has a nonconstant section, a contradiction.
\par\noindent(ii) $y_0y_5y_{10}$ is a section of $\det W_E=\cal O_E(45p_0)$
which is invariant under the induced action of $\Delta$ on $\det W_E$.
It defines an invariant divisor on $E$ whose image in the quotient is a
divisor linearly equivalent to $5p_0$.
\end{Proof}

Let us now look at the corresponding action of $\Bbb Z_3\times\Bbb Z_3$
on the trivial projective bundle $E\times\Bbb P^2=E\times\Bbb P(V\spcheck )$.
Its quotient is a $\Bbb P^2$-bundle
\[
  \Bbb P^2_E=\Bbb P(\cal E_E)
\]
where we use the geometric projective bundle. By the above lemma
$\Bbb P^2_E$ is the unique indecomposable $\Bbb P^1$-bundle
over $E$ with invariant $e=-1$ \cite[V. theorem 2.15]{Ha}.
We consider the quotient map
\[
  \pi\colon E\times\Bbb P^2\to\Bbb P_E^2.
\]

Clearly this map is unramified and we can use $\pi$ to compute the
cohomology of line bundles on $\Bbb P_E^2$. This was done in \cite{CC}
for the dual bundle $\cal E_E\spcheck$. We are particularly interested
in line bundles numerically equivalent to the anticanonical bundle.
The Picard group $\Bbb P_E^2$ is generated by the tautological bundle
$\cal O_{\Bbb P_E^2}(1)$ and the pullback of the Picard group on $E$. The
pullback of any line bundle on  $\Bbb P_E^2$ to $E\times\Bbb P^2$
is the tensor product of a line bundle on $E$ and a line bundle
on $\Bbb P^2$.

\begin{Lemma}[Catanese, Ciliberto]\label{L3}
  If $\cal O_{\Bbb P_E^2}(L)$ is numerically equivalent to the anticanonical
  bundle $\cal O_{\Bbb P_E^2}(-K)$ and $h^\circ(\cal O_{\Bbb P^2_E}(L))>0$
  then either $L\equiv-K$ in which case
  $h^\circ(\cal O_{\Bbb P_E^2}(L))=2$
  or $3L\equiv-3K$ and $L\not\equiv-K$ in which case
  $h^\circ(\cal O_{\Bbb P_E^2}(L))=1$. Moreover there are 8 nonisomorphic
  bundles of the latter kind corresponding to the nontrivial characters
  of $\Bbb Z_3\times\Bbb Z_3$.
\end{Lemma}

\begin{Proof}
Let $L=-K+\rho$ where $\rho$ is the pullback of a degree 0 line bundle on $E$.
Since $\pi$ is unramified $\pi^*(-K+\rho)=-K_{E\times\Bbb P_2}+\rho'$
where $\rho'$ also has degree 0. This bundle can only have sections when
$\rho'=0$. Since
\[
  \cal O_{\Bbb P_E^2}(-K+\rho)\subset\pi_*\pi^*\cal O_{\Bbb P_E^2}(-K+\rho)
  =\pi_*(-K_{E\times\Bbb P_2})
\]
we are left to consider the decomposition
\[
  \pi_*(-K_{E\times\Bbb P_2})
  =\bigoplus_{\Chi\in(\Bbb Z_3\times\Bbb Z_3)\spcheck}(-K+L_\Chi)
\]
where $L_\Chi$ is the torsion bundle associated to the character $\Chi$.
I.e., $L_\Chi$ is a torsion bundle of degree 0 on $E$. Hence $L$ is of
the form stated. The sections of $L$ are given by the sections of
$\cal O_{E\times\Bbb P^2}(-K_{E\times\Bbb P_2})$ associated to the
character $\Chi$. By what we have said in the previous paragraph the
dimension of these sections is 2 if $\Chi$ is trivial and 1 otherwise.
\end{Proof}

\begin{Lemma}\label{L4}
  In the pencil $|-K|$ the singular members are four singular scrolls
  while the smooth members are abelian surfaces among which one is
  isomorphic to $E\times E$. The divisors $-K+L_\Chi$, $\Chi$ nontrivial,
  are smooth bielliptic surfaces.
\end{Lemma}

\begin{Proof}
The divisors $\pi^*(-K)$ and $\pi^*(-K+L_\Chi)$ are $E\times E_{\lambda'}$
where $E_{\lambda'}$ is a member of the Hesse pencil and
$E\times F_{(a,b)}$, respectively. On each of these surfaces the
$\Bbb Z_3\times\Bbb Z_3$-action is the one described in the previous paragraph.
When $E_{\lambda'}$ is smooth then the 9 base points of the Hesse
pencil form a subgroup of the product, so the quotient is abelian.
In particular when $\lambda'=\lambda$ we get $E\times E/\Bbb Z_3\times\Bbb Z_3$
where $\Bbb Z_3\times\Bbb Z_3$ acts diagonally. It is easy to see that this
quotient is again isomorphic to $E\times E$ (we shall soon discuss this
in more detail). When $E_{\lambda'}$ is a triangle then the surface
upstairs is the union of three scrolls, whose quotient downstairs is
irreducible since the group acts transitively on the edges of the triangles.
Finally on $E\times F_{(a,b)}$ the group acts with translation on the first
factor and with translation and multiplication on the second factor.
Hence the quotient is bielliptic.
\end{Proof}

We want to describe the intersection of the abelian and bielliptic surfaces
with the special abelian surface $A_0\iso E\times E$ described above.
By $A_K$ we'll denote the general abelian surface in  $|-K|$. Let $T_K$
be the singular scrolls in $|-K|$ and $B_{(a,b)}$ the bielliptic
surfaces. Let us first consider the structure of the abelian and bielliptic
surfaces. Each of them has an elliptic fibration over $E$ whose fibres
are the plane cubic curves $E_{\lambda'}$, $F_{(a,b)}$ respectively.
For the abelian surfaces the elliptic fibration over $E_{\lambda'}$
upstairs remains an elliptic fibration over
$E_{\lambda'}/\Bbb Z_3\times\Bbb Z_3=E_{\lambda'}$ downstairs, the fibres
being isomorphic to $E$. Upstairs the intersection of $E\times E$ and
$E\times E_{\lambda'}$ is 9 translates of the curve $E$ over the 9 base
points of the Hesse pencil. On the quotient these translates are mapped
to the same curve isomorphic to $E$, which in turn is a member of the
fibration over $E_{\lambda'}$ described above. For the bielliptic surfaces
the elliptic fibration over $F_{(a,b)}$ upstairs is mapped to an
elliptic fibration over $F_{(a,b)}/\Bbb Z_3\times\Bbb Z_3\iso\Bbb P^1$
downstairs, i.e., a pencil. The intersection upstairs with $E\times E$
is 9 translates by the group of $E$ over the 9 points of intersection
$E_\lambda\cap F_{(a,b)}$ which are mapped to the same curve
$A_0\cap B_{(a,b)}$ downstairs. This is a member of the pencil over
$\Bbb P^1$. The intersection $A_K\cap B_{(a,b)}$ downstairs is linearly
equivalent to and different
from $A_0\cap B_{(a,b)}$, so $|-K|$ restricts to $B_{(a,b)}$ to give
the pencil described above. In particular the three
triple fibres of this pencil are the intersections $T_K\cap B_{(a,b)}$
for the scrolls $T_K$ coming from triangles $T_{(i,j)}$ with
$(i,j)\ne\pm(a,b)$.

Our next aim is to describe the intersection of $A_K$, resp.\ of
$B_{(a,b)}$ with $A_0$ more arithmetically. We look at the map
\[
  \begin{pmatrix}\phantom{-}3&\phantom{-}0\\
                           -2&          -1
  \end{pmatrix}\colon
  \left\{\begin{aligned}E\times E &\to     E\times E\\
                        (q_1,q_2) &\mapsto (3q_1,-2q_1-q_2).
  \end{aligned}\right.
\]
The kernel of this map is the group
$E^{(3)}$ of $3$-torsion points of $E$ embedded diagonally into $E\times E$.
Hence the above map induces an isomorphism
\[
  A_0=(E\times E)/\Bbb Z_3\times\Bbb Z_3\iso E\times E.
\]
Whenever we shall refer to $A_0$ as a product it will be via this
isomorphism. Note that the curve $\{(q,-2q);\ q\in E\}$ goes
9$:$1 onto the first factor and that $\{(0,-q);\ q\in E\}$ is mapped
isomorphically onto the second factor. Moreover the curve
$\{(q,-5q);\ q \in E\}$ goes 9$:$1 onto the diagonal and
$\{(q,q);\ q \in E\}$ is mapped 9$:$1 onto the antidiagonal of $E\times E$.
Finally we consider the map given by $\left(\begin{smallmatrix}-2&-1\\
-5&2\end{smallmatrix}\right)$ upstairs. One checks immediately that this
induces an endomorphism downstairs, and that this endomorphism is
$\left(\begin{smallmatrix}0&3\\3&0\end{smallmatrix}\right)$, i.e., 3 times
the standard involution interchanging the factors of $E\times E$.

\begin{Lemma}\label{L5}
Let $\Delta_E$ be the diagonal in $A_0=E\times E$.

\noindent\rom{(i)} The curve $A_K\cap E\times E$ is
  \[
    \{(q,r);\ 3r+2q=0\}
  \]
and  $A_K\cap \Delta_E$ consists of the 25 points
  \[
     \{(p,p);\ 5p=0\}.
  \]

\noindent\rom{(ii)} The curves $B_{(a,b)}\cap E\times E$ are
  \[
    \{(q,r);\ 3r+2q=-\tau_{(a,b)}\}
  \]
and $B_{(a,b)}\cap \Delta_E$ are the sets of points
  \[
     \{(p,p);\ 5p=-\tau_{(a,b)}\},
  \]
where $0\ne\tau_{(a,b)}$, $3\tau_{(a,b)}=0$.

\end{Lemma}

\begin{Proof}
(i) Upstairs $E\times E_\lambda\cap E\times E=\{(q,\tau_3);\ q\in E,
\ 3\tau_3=0\}$. The image of this set downstairs is $\{(3q,-2q-\tau_3);
\ q\in E\}$ which is the curve described. The second part follows
immediately.

\noindent (ii) Similarly $E\times F_{(a,b)}\cap E\times E=\{(q,\tau_9);
\ q\in  E,\ 3\tau_9=\tau_{(a,b)}\}$ where the $\tau_{(a,b)}$ are the
3-torsion points on $E$. Downstairs this is
$\{(3q,-2q-\tau_9);\ q\in E\}$ which gives the claim.
\end{Proof}

At this point we want to return to the product $E\times\Bbb P^2$. Let
$p,q$ be the projections onto $E$ and $\Bbb P^2$. Let
$\cal O_E(15p_0)\boxtimes\cal O_{\Bbb P^2}(1)=p^*\cal O_E(15p_0)\otimes
q^*\cal O_{\Bbb P^2}(1)$. The centre of the diagonal
$\Delta\subset H_3\times H_3$ acts trivially on this line bundle which,
therefore, descends to a line bundle $\cal L$ on $\Bbb P_E^2$.

\begin{Proposition}\label{P6}
  \rom{(i)} $h^\circ(\cal L)=5$\par
  \noindent\rom{(ii)} The following sections are invariant under
  $\Delta$, hence define a basis of $\mathrm H^\circ(\cal L)$:
  \begin{align*}
    s_0 &= y_{ 0}\otimes x_0+y_{ 5}\otimes x_1+y_{10}\otimes x_2\\
    s_1 &= y_{ 3}\otimes x_0+y_{ 8}\otimes x_1+y_{13}\otimes x_2\\
    s_2 &= y_{ 6}\otimes x_0+y_{11}\otimes x_1+y_{ 1}\otimes x_2\\
    s_3 &= y_{ 9}\otimes x_0+y_{14}\otimes x_1+y_{ 4}\otimes x_2\\
    s_4 &= y_{12}\otimes x_0+y_{ 2}\otimes x_1+y_{ 7}\otimes x_2
  \end{align*}
\end{Proposition}

\begin{Proof}
(i) Clearly
\[
  \mathrm H^\circ(\cal O_E(15p_0)\boxtimes\cal O_{\Bbb P^2}(1))
  = \mathrm H^\circ(\cal O_E(15p_0))\otimes
    \mathrm H^\circ(\cal O_{\Bbb P^2}(1)).
\]
As an $H_3$-module
\[
  \mathrm H^\circ(\cal O_E(15p_0)) = 5 V\spcheck.
\]
This can be seen by looking at the subspaces spanned by
$(y_0,y_5,y_{10})$, $(y_3,y_8,y_{13})$, $(y_6,y_{11},y_1)$,
$(y_9,y_{14},y_4)$, $(y_{12},y_2,y_7)$. Hence as an
$H_3$-module
\[
  \mathrm H^\circ(\cal O_E(15p_0)\boxtimes\cal O_{\Bbb P^2}(1))
  = 5V\spcheck\otimes V
  = 5(\bigoplus_{\Chi\in(\Bbb Z_3\times\Bbb Z_3)\spcheck}V_\Chi).
\]
I.e. we have 5 invariant sections, and thus $h^\circ(\cal L)=5$.

\noindent(ii) It is straightforward to check that the $s_i$ are invariant
under $H_3$.
\end{Proof}

Next we consider the subgroup of $H_{15}$ spanned by $\sigma_{15}^3,
\tau_{15}^3$. From \eqref{(3)}
\[
  [\sigma_{15}^3,\tau_{15}^3]
  = \varepsilon_{15}^{-9}\cdot\operatorname{id}
  = \varepsilon_5^{-2}\cdot\operatorname{id}
  \quad (\varepsilon_5=e^{2\pi i/5}).
\]
Hence mapping $\sigma_{15}^3$ to  $\sigma_5$ and $\tau_{15}^3$ to
$\tau_5$ we can identify this subgroup with
$H_5\subset\operatorname{SL}(\Bbb C^5)$, where this inclusion is given
by the representation which arises from the Schr\"odinger representation
of the Heisenberg group $H_5$ of level 5 by replacing $\varepsilon$
by $\varepsilon^2$. Now let $H_5$ act on
$E\times\Bbb P^2$ where the action on the second factor is trivial.
Then $H_5$ acts on $\cal O_E(15p_0)\boxtimes\cal O_{\Bbb P^2}(1)$
and it is straightforward to check that this action commutes with $H_3$.
Hence we get an action of $H_5$ on $\cal L$.

\begin{Proposition}\label{P7}
  The action of $H_5$ on $\mathrm H^\circ(\cal L)$ is given by
  \[
    \sigma_5(s_i)=s_{i-1},\quad\tau_5(s_i)=\varepsilon_5^{-2i}s_i.
  \]
\end{Proposition}

\begin{Proof}
Straightforward calculation.
\end{Proof}

We have involutions on  $E$ (given by $x\mapsto-x$) and on $\Bbb P^2$
(given by $\iota_3(x_i)=x_{-i}$). Hence we have an involution $\iota$ on
$E\times\Bbb P^2$. This lifts to an involution on
$\cal O_E(15p_0)\boxtimes\cal O_{\Bbb P^2}(1)$ where it acts on sections by
\begin{equation}\label{(4)}
  \iota(y_i\otimes x_i) = y_{-i}\otimes x_{-i}.
\end{equation}
This involution does not commute with $H_3$, but we have an action
of a semi-direct product $H_3\rtimes\langle\iota\rangle$.
In the quotient this defines an involution on $\Bbb P_E^2$ and on $\cal L$.
Note that on $A_0=E\times E$ this is given by
$\left(\begin{smallmatrix}-1&\phantom{-}0\\
\phantom{-}0&-1\end{smallmatrix}\right)$.

\begin{Proposition}\label{P8}
  $\iota$ acts on $\mathrm H^\circ(\cal L)$ by
  \[
    \iota(s_i)=s_{-i}.
  \]
\end{Proposition}

\begin{Proof}
Immediately from \eqref{(4)}.
\end{Proof}

Finally we remark that we really have an action of $(\Bbb Z_3\times\Bbb Z_3)^2$
on $E\times\Bbb P^2$ and that $\Bbb P_E^2$ was constructed by taking the
quotient
with respect to the diagonal. Hence we have still got an action of
$\Bbb Z_3\times\Bbb Z_3$ on $\Bbb P_E^2$ which on every fibre of $\Bbb P_E^2$
lifts to the Schr\"odinger representation of $H_3$.

Let $\Delta_E$ be the diagonal in $A_0=E\times E$. We can consider the
blow-up
\[
  \rho\colon U\to\Bbb P_E^2
\]
along $\Delta_E$. Since $\Delta_E$ is a section of $\Bbb P_E^2$ the
variety $U$ has the structure of an $\Sigma^1$-bundle over $E$. Here
$\Sigma^1$ denotes the $\Bbb P^1$-bundle over $\Bbb P^1$ with $e=-1$.

\begin{Lemma}\label{L9}
  $U$ has the structure of a $\Bbb P^1$-bundle over $S^2E$.
\end{Lemma}

\begin{Proof}
Let $E_\Delta$ be the exceptional surface over $\Delta_E$
and $B=\cal O_{\Bbb P_E^2}(1)$. By $F$ we denote the class of a fibre
of $\Bbb P_E^2$. For $\beta$ sufficiently large $|B-E_\Delta+\beta F|$
is base point free. This linear system maps each $\Sigma^1$
to a $\Bbb P^1$-bundle
over a scroll over $E$, and it remains to determine this scroll. To do this
we look at $A_0=E\times E$. The map given by $|B-E_\Delta+\beta F|$
restricted to a curve $\{q\}\times E$ is nothing but projection of the
plane cubic $E\subset\Bbb P^2$ from the point $q$. Hence we
get an involution on $E\times E$ whose branch locus is the curve
$\Delta'=\{(q,t)\in E\times E;\ 2t+q=0\}$. $\Delta'$ is the
image of $E\to E\times E$, $q\mapsto\left(\begin{smallmatrix}2q\\-q
\end{smallmatrix}\right)$. The isomorphism of $E\times E$ given by the
matrix $\left(\begin{smallmatrix}-1&-1\\
\phantom{-}0&\phantom{-}1\end{smallmatrix}\right)$
maps $\Delta'$ to the diagonal and the curves $\{q\}\times E$ to the
translates of the antidiagonal $\{(q,-q);\ q\in E\}$. Under this
isomorphism the above involution becomes the standard involution given
by interchanging the factors. This proves that the scroll in question
is indeed $S^2E$. We have indeed a (locally trivial) $\Bbb P^1$-bundle
by \cite[V.4.1]{BPV} and the remark after this.
\end{Proof}

In view of the above lemma we shall change our notation and write
\[
  \Bbb P_{S^2E}^1=U.
\]

The strict transform of $A_0=E\times E$ under the map
$\rho\colon\Bbb P_{S^2E}^1\to\Bbb P_E^2$ is again $E\times E$. The other
abelian surfaces $A_K$ in $\Bbb P_E^2$ are blown up in the points
\[
  \Delta_E\cap(A_K\cap A_0) = \{q;\ 5q=0\},
\]
i.e., in 25 distinct points (cf. lemma \ref{L5}). Similarly the surfaces
$B_{(a,b)}$ are blown up in 25 distinct points by $\rho$.

\section{Two quintic hypersurfaces in $\Bbb P^4$}\label{Par2}

Let $E$ be an elliptic normal curve of degree 5 in $\Bbb P^4$, embedded by
the linear system $|5p_0|$, where $p_0 \in E$ is the origin which we have
chosen before. We assume that $E$ is invariant under the action of
$H_5$ given by its Schr\"odinger representation.
In this paragraph we will describe geometrically the following diagram
\unitlength1.2pt
\[
  \begin{picture}(70,50)
  \put(0,40){$\Bbb P^1_{S^2E}$}  \put(32,43){$^\rho$}  \put(60,40){$\Bbb
P_E^2$}
                        \put(12,43){\vector(1,0){46}}

       \put(20,31){$^{f_1}$}
      \put(12,36){\vector(2,-1){20}}

                            \put(33,20){$\cal F$}
  \put(-1,21){$^{\pi_1}$}                           \put(67,21){$^{\pi_2}$}
  \put(6,36){\vector(0,-1){28}}                    \put(66,36){\vector(0,
-1){28}}

        \put(18,13.5){$^{p_1}$}             \put(50,13){$^{p_2}$}
       \put(31,18){\vector(-2,-1){20}}    \put(42,18){\vector(2,-1){20}}

  \put(2,0){$V_1$}                                      \put(62,0){$V_2$.}
  \end{picture}
\]
Here $V_1$ and $V_2$ are quintic hypersurfaces in $\Bbb P^4$: $V_1$ is the
secant variety to $E$ and $V_2$ is ruled in an elliptic family of planes.
$V_2$ is singular (set theoretically) along an elliptic quintic scroll
whose trisecant variety it is. $\Bbb P^1_{S^2E}$,  $\Bbb P_E^2$ and $\rho$
are as in
\ref{Par1}. Via $\pi_1$ and $\pi_2$ they are minimal desingularizations
of $V_1$ and $V_2$. The varieties above are tied together via an incidence
variety $\cal F$ consisting of pairs $(p,Q)$ where $Q$ is a singular quadric
containing $E$ and $p\in\operatorname{Sing}Q$. The morphisms to $V_1$ and
$V_2$ are projections to the first and second factor respectively.
Hence both hypersurfaces can be described in terms of the 3-dimensional
family of singular quadrics through $E$. Furthermore one can make the
identifications
$$\PP^1_{S^2E} = \{(p,\{e_1,e_2\}) | e_1, e_2\in E, p\in L_{\langle e_1,
e_2\rangle}\}
=\{(p,W_p) | p\in V_1, W_p=\mathrm H^\circ (I_{E\cup \{p\}}(2))\}$$

where $L_{\langle e_1,e_2\rangle}$ is the secant line through $e_1,e_2$.
The second projection for these incidence varieties is the map to
$S^2E\subset{\Bbb P(\mathrm H^\circ(I_E(2))^{\spcheck})}$. Thus $\PP^1_{S^2E}$
is the graph of the quadro cubic Cremona transformation restricted to the
secant
variety $V_1$ of $E$.  One can also show that $V_2\subset
{\Bbb P(\mathrm H^\circ(I_E(2)))}$ in this setting is the natural dual
to $S^2E\subset{\Bbb P(\mathrm H^\circ(I_E(2))^{\spcheck})}$.
 Since this will not be essential for our argument we will omit the proofs.
Some
of the results we collect here are also contained in \cite{EL}, \cite{Hu},
\cite{d'Al}.

Let $E$ be as above. It is well known that $h^\circ(I_E(2))=5$ and that a
basis for the space of quadrics through $E$ is given by
\[
  Q_i = x_i^2+ax_{i+2}x_{i+3}-\frac1ax_{i+1}x_{i+4}\quad(i\in\Bbb Z_5).
\]
Here $a\in\Bbb C\cup\{\infty\}$ and five such quadrics define a smooth
elliptic curve if and only if $a$ is not a vertex of the icosahedron,
i.e., $a\ne0,\infty,\varepsilon_5^k(\varepsilon_5^2+\varepsilon_5^3),
\varepsilon_5^k(\varepsilon_5+\varepsilon_5^4)$, $k=0,\ldots,4$
($\varepsilon_5=e^{2\pi i/5}$).

\begin{Definition}\label{D10}
  \rom{(i)} For $y\in\Bbb P^4$ let $M(y)$ be the symmetric 5$\times$5-matrix
  \[
    M(y):=(y_{i+j}z_{i-j})\quad0\le i,j\le4
  \]
  where $z\in\Bbb P^4$, $z_i=z_{-i}$ and $z_0=2$, $z_1=a$, $z_2=-\frac1a$.

  \noindent\rom{(ii)} Let $\cal F$ denote the incidence variety
  \[
    \cal F:=\{(x,y)\in\Bbb P^4\times\Bbb P^4;\ M(y)\tr x=0\}
  \]
  and let $V_1$ and $V_2$ denote the images under the first and second
projections
  $p_1$ and $p_2$ of $\cal F$ to the respective $\Bbb P^4$'s.

  \noindent\rom{(iii)} Let $M'(x)$ be the 5$\times$5-matrix defined by
  \[
    M'(x)\tr y=M(y)\tr x.
  \]
\end{Definition}

\begin{Remark}
  This set-up was also considered in \cite{A} in the case of a general point
  $z\in\Bbb P^4$. Here we have chosen a special point, namely one that lies
  on the conic section invariant under the icosahedral group $A_5$ on the Bring
  plane $z_i=z_{-i}$; $i=0,\ldots,4$ (see \cite{BHM}). This conic can be
  identified naturally with the modular curve of level 5, $X(5)$, which is in
  1$:$1-correspondence with $H_5$-invariantly embedded elliptic quintics
  in $\Bbb P^4$. Under this identification $z$ corresponds to the curve $E$
  we have started with. The matrix $M$ was first considered by Moore.
\end{Remark}

\begin{Proposition}\label{P12}
  \rom{(i)} The set of quadrics $\{xM(y)\tr x;\ y\in\Bbb P^4\}$ is
  $\Bbb P(\mathrm H^\circ(I_E(2)))$.

  \noindent\rom{(ii)} $\cal F$ can be identified with the incidence variety
  of pairs $(p,Q)$ where $Q$ is a singular quadric through $E$ and
  $p\in\operatorname{Sing}Q$.

  \noindent\rom{(iii)} $V_1$ and $V_2$ are Heisenberg invariant quintic
  hypersurfaces in $\Bbb P^4$.
\end{Proposition}

\begin{Proof}
(i) It is easily checked that
\[
  2Q_{3i}(x)=xM(e_i)\tr x\quad(i=0,\ldots,4).
\]

\noindent(ii) This follows immediately since $\cal F$ is given by
$M(y)\tr x=0$.

\noindent(iii) Since $M(y)\tr x=0$ is equivalent to $M'(x)\tr y=0$
it follows that $V_1$ and $V_2$ are given by the quintic equations
$\operatorname{det}M'(x)=0$ resp.\ $\operatorname{det}M(y)=0$. Because both
$M'(x)$ and $M(y)$ are invariant under $H_5$, up to an even number of
permutations of rows and columns, these equations are $H_5$-invariant.
Since $H_5$ has no characters on $\mathrm H^\circ(\cal O_{\Bbb P_4}(n))$,
$n<5$ it follows that $V_1$ and $V_2$ are in fact reduced of degree 5.
\end{Proof}

\begin{Remark}\label{R13}
  Since the general singular quadric through $E$ has rank 4, it follows that
  the projection $\cal F\to V_2$ is generically finite and hence $\cal F$
  is also of dimension 3.
\end{Remark}

\begin{Corollary}\label{C14}
  \rom{(i)} $V_1$ is the locus of singular points of the singular quadrics
  through  $E$.

  \noindent\rom{(ii)} $V_1=\operatorname{Sec}E$.
\end{Corollary}

\begin{Proof}
(i) is obvious.

\noindent\rom{(ii)} If $p\in\operatorname{Sec}E\backslash E$ then projection
to  $\Bbb P^3$ from $p$
maps $E$ to a nodal quintic curve in $\Bbb P^3$, which always lies on a
quadric surface. Hence $p$ lies on a quadric cone through $E$. This
implies $\operatorname{Sec}E\subset V_1$ and since both hypersurfaces have
degree 5, the claim follows.
\end{Proof}

\begin{Corollary}\label{C15}
$V_2$ is the discriminant locus of the family of quadrics through  $E$.
\end{Corollary}

\begin{Proof}
Clear.
\end{Proof}

The mapping $p_1$ (resp.\ $p_2$) is a ``small resolution'' of $V_1$
(resp.\ $V_2$), i.e., a (singular) point where $M'(x)$ (resp.\ $M(y)$)
has rank 3 is replaced by a $\Bbb P^1$. One of our aims is to describe
the rank 3 loci. For $M(y)$ the corresponding $\Bbb P^1$s yield the locus
of the singular lines. We shall come back to this later. For $M'(x)$ it is
simpler.

\begin{Proposition}\label{P16}
  \rom{(i)}
  The quintic hypersurface $V_1 = \operatorname{Sec} E$ is singular
  precisely at $E$ where it has multiplicity 3.

  \noindent\rom{(ii)}
  The curve $E$ is exactly the locus where rank $M'(x) = 3$.
\end{Proposition}

\begin{Proof}
  We shall first prove that the multiplicity of $\operatorname{Sec} E$
  along $E$ is three. This was already known to Segre \cite{Seg}, \cite{Sem}.
  Here we reproduce his proof. We consider a point $p \in E$ and
  choose a general line $l$ through $p$. We can assume that $l$ meets
  $\operatorname{Sec} E$ transversally at a finite number of smooth points
  outside $p$. Since secants and tangents of $E$ do not meet outside $E$
  (see \cite[Lemma IV.11]{Hu}) every such point of intersection lies on a
  unique secant or tangent of $E$. On the other hand projection from a
  general line $l$ maps $E$ to a plane curve of degree 4 which, by the
  genus formula, must have 2 nodes. Hence $l$ intersects
  $\operatorname{Sec} E$ in precisely 2 points (counted properly) outside
  $p$, and it follows that the multiplicity of $\operatorname{Sec} E$
  along $E$ is 3.

  Now assume that a singularity $x$ of $\operatorname{Sec} E$ exists
  outside $E$. Let $l$ be a line through $x$ which meets $E$ in a point
  $p$, but is neither a secant nor a tangent line of $E$. (Such a line
  exists since $x$ lies on at most one secant or tangent --- see above).
  Using \cite[Proposition IV.4.6]{Hu} we can also assume that $l$ is not a
  singular line of a rank 3 quadric through $E$. It follows that $l$ is not
  contained in $\operatorname{Sec} E$. The latter would only be possible if
  projection from $l$ defines a 2:1 map onto a conic, but this implies that
  $l$ is the vertex of a rank 3 quadric. Since the multiplicity of
  $\operatorname{Sec} E$ along $E$ is 3 and since $x$ was assumed to be
  singular, it follows that the intersection of $l$ and
  $\operatorname{Sec} E$ consists precisely of the two points $x$ and $p$.
  Projection from $l$ now gives a curve of degree 4 and genus 1 in $\PP^2$
  with exactly one singular point, given by the unique secant or tangent of
  $E$ through $x$. On the other hand, we can project from $x$ first. In
  this case $E$ is mapped to a quintic curve $E^{\prime}$ in ${\PP^3}$
  with one singularity and arithmetic genus 2, which lies on a unique
  quadric surface $Q'\subset\PP^3$. We have two possibilities

  \noindent\rom{(1)}
  $Q'$ is a smooth quadric. In this case $E'$ is a divisor on $Q'$ of
  bidegree $(2,3)$. Then projection from a general point on $E'$ (which
  corresponds to a general choice of the point $p \in E$) projects $E'$
  to a plane quartic with 2 different singularities, a contradiction to
  what we have found above.

  \noindent\rom{(2)}
  $Q'$ is a quadric cone. In this case $E'$ contains the vertex of this
  cone as a smooth point and meets every ruling of $Q'$ in 2 points
  outside the vertex. Again projection from a general point on $E'$
  gives a quartic plane curve with two different singularities and we
  have arrived at the same contradiction as above.

  It follows that $\operatorname{Sec} E$ has no singularities outside $E$.

  \noindent\rom{(ii)}
  The locus where rank $M'(x)$ has rank $\le 3$ is contained in
  $\operatorname{Sing} V_1$. On the other hand if $p \in E$, then
  projection from $p$ gives a smooth quartic elliptic curve in $\PP^3$
  which lies on a pencil of quadrics. Hence $E$ lies precisely on a
  pencil of quadric cones with vertex $p$. It follows that
  \[
    M(y) \tr{p} = M'(p) \tr{y}
  \]
  for $y$ in some (linear) $\PP^1$, and $M'(p)$ has rank exactly 3.
\end{Proof}

\begin{Remark}\label{R17}
One shows easily that
\[
  M'(x)=\fracwithdelims(){\partial Q_{3j}}{\partial  x_i}_{0\le i,j\le4}.
\]
\end{Remark}

\medskip

We consider the natural desingularization
\[
  \widetilde V_1 := \{(p,\{e_1,e_2\})\in\Bbb P^4\times S^2E;
  \ p\in L_{\langle e_1,e_2\rangle}\}
\]
of $V_1=\operatorname{Sec}(E)$, where $L_{\langle e_1,e_2\rangle}$ is
the secant line through $e_1,e_2$.
Projection onto $S^2E$ gives $\widetilde V_1$ a structure of a $\Bbb
P^1$-bundle
over
the surface $S^2E$. Let $\pi_1$ denote projection to the first factor.
Then $\pi_1$ contracts the divisor
\[
  D_1:=\{(p,\{e_1,e_2\});\ p\in\{e_1,e_2\}\}.
\]
We have an isomorphism
\[
  \psi_1\colon
  \left\{\begin{aligned} D_1             &\iso  E\times E\\
                         (p,\{e_1,e_2\}) &\mapsto (p,e_1+e_2).
  \end{aligned}\right.
\]
Next we consider the natural composition
\[
  \widetilde V_1\to S^2E\to E
\]
where the map $S^2E\to E$ maps $\{e_1,e_2\}$ to $e_1+e_2$. The fibre of
this map over a point $e\in E$ is the surface
\[
  \{(p,\{e_1,e_2\});\ e_1+e_2=e,\ p\in L_{\langle e_1,e_2\rangle}\}.
\]
This is a ruled surface over the curve
\[
  E/\kappa\iso\Bbb P^1
\]
where $\kappa$ is the involution on $E$ given by $\kappa(q)=-q+e$.
Via $\pi_1$ this
is a smooth, rational ruled surface in $\Bbb P^4$, i.e., a cubic scroll.
As an abstract surface this is $\Bbb P^2$ blown up in a point, or
equivalently the Hirzebruch surface $\Sigma^1$. In this way $\widetilde V_1$
acquires the structure of a $\Sigma^1$-fibration over  $E$. We denote the
fibre of this fibration over a point $e\in E$ by $\Sigma_e^1$. We shall
often identify $\Sigma_e^1$ with $\pi_1(\Sigma_e^1)$. Thus we can write
\[
\widetilde V_1=\{(p,e);\ p\in \Sigma_e^1\subset V_1\}.
\]
 We will use this notation in the sequel.

\begin{Proposition}\label{P18}
  The map $\pi_1$ defines an isomorphism from $\widetilde V_1\sm D_1$ with
  $\operatorname{Sec} E \sm E$. It contracts $D_1$ to the curve $E$ and
  its differential has rank 2 at every point of $D_1$.
\end{Proposition}

\begin{Proof}
  Since every point on $\operatorname{Sec} E \sm E$ lies on a unique
  secant or tangent of $E$ the map from $\widetilde V_1 \sm D_1$ to
  $\operatorname{Sec} E \sm E$ is bijective. Since both are smooth, it
  is an isomorphism. We have already seen that $\pi_1$ contracts $D_1$
  to the curve $E$. Hence the differential of $\pi_1$ along $D_1$ has
  rank at most 2. On the other hand consider the fibres $\Sigma_e^1$ of
  the map $\widetilde V_1 \to E$. Via the map $\pi_1$ they are embedded into
  $\PP^4$, and hence the differential of $\pi_1$ has rank at least 2 at
  every point of $\widetilde V_1$.
\end{Proof}

Now we return to the cubic scroll $\Sigma_e^1\subset\Bbb P^4$. Since
$\Sigma_e^1$ is the degeneration locus of a 2$\times$3 matrix with linear
coefficients, it follows that there is a $\Bbb P^2$ of quadrics
containing $\Sigma_e^1$. All of these quadrics are singular. Geometrically
they arise as follows: Projection from $p\in\Sigma_e^1$ maps $\Sigma_e^1$
to  a quadric in $\Bbb P^3$. Then take the cone over the quadric in
$\Bbb P^3$. Note that the quadric surface in $\Bbb P^3$ is singular if
and only if $p$ is on the exceptional line in $\Sigma_e^1$.
In this case the corresponding quadric hypersurface is singular
along the exceptional line in  $\Sigma_e^1$. Finally it follows easily from
$V_1=\operatorname{Sec}E$ that every singular quadric through  $E$
arises in the way described above.

We define
\[
  \widetilde V_2 := \{(Q_e,e);\ e\in E,
  \ \text{$Q_e$ is a quadric through $\Sigma_e^1$}\}.
\]
Via the obvious map $\widetilde V_2\to E$ this carries the structure of a
$\Bbb P^2$-bundle. For $p\in\Sigma_e^1$ we denote by $Q_e=Q_e(p)$ the
unique quadric through $\Sigma_e^1$ which is singular at $p$.
This enables us to define the following maps:
\[
  f_1\colon\left\{\begin{aligned}\widetilde V_1 &\to     \cal F\\
                                 (p,e)      &\mapsto (p,Q_e)
           \end{aligned}\right.
\]
where $(p,e)$ stands for the point $p\in\Sigma_e^1$, and
\[
  f_2\colon\left\{\begin{aligned}\widetilde V_1 &\to     \widetilde V_2\\
                                 (p,e)      &\mapsto (Q_e,e).
           \end{aligned}\right.
\]
In this way we get a commutative diagram
\unitlength1.2pt
\begin{equation}\label{(D)}
  \begin{picture}(70,50)
  \put(2,40){$\widetilde V_1$}  \put(32,43){$^{f_2}$}  \put(62,40){$\widetilde
V_2$}
                        \put(12,43){\vector(1,0){46}}
       \put(20,31){$^{f_1}$}
      \put(12,36){\vector(2,-1){20}}
                            \put(33,20){$\cal F$}
  \put(-1,21){$^{\pi_1}$}                           \put(67,21){$^{\pi_2}$}
  \put(6,36){\vector(0,-1){28}}                    \put(66,36){\vector(0,
-1){28}}
        \put(18,13.5){$^{p_1}$}             \put(50,13){$^{p_2}$}
       \put(31,18){\vector(-2,-1){20}}    \put(42,18){\vector(2,-1){20}}
  \put(2,0){$V_1$}                                      \put(62,0){$V_2$.}
  \end{picture}
\end{equation}
Moreover it follows from our geometric discussion
that $f_2$ contracts precisely the divisor
\[
  X:=\{(p,e);\ p\in\ \text{exceptional line in $\Sigma_e^1$}\}.
\]
In other words $f_2$ is the blowing down map from the $\Sigma^1$-bundle
$\widetilde V_1$ to the $\Bbb P^2$-bundle $\widetilde V_2$. Furthermore
$X$ is an elliptic ruled surface and $\pi_1(X)$ is the locus of singular
lines.

We now return to the divisor $D_1\iso E\times E$ in $\widetilde V_1$.

\begin{Lemma}\label{L18}
\rom{(i)} $f_1$ is an isomorphism outside $D_1$.

\noindent\rom{(ii)} $f_1(p,e)=f_1(p',e')$ if and only if
$p'=p\in E$ and $e+e'=-p$.
\end{Lemma}

\begin{Proof}
\rom{(i)}
This follows from proposition \ref{P16} \rom{(i)}.

\noindent\rom{(ii)}
If $f_1(p,e) = f_1(p',e')$ then clearly $p= p'$ by construction of the map
$f_1$.
Now $Q_e = Q_{e'}$ means that $Q_e$ is a singular quadric with vertex $p\in E$
containing both $\Sigma_e^1$ and $\Sigma_{e'}^1$.
$\Sigma_e^1$ and $\Sigma_{e'}^1$ are determined by the families of planes
in $Q_e$: A plane intersects $E$ in two points besides $p$, defining
a line in the ruling of the scroll. If $L_{\langle e_1,e_2\rangle}\subset
\Sigma_e^1$ and $L_{\langle e_1',e_2'\rangle}\subset\Sigma_{e'}^1$ then
$e_1,e_2,e_1',e_2'$ and $p$ are contained in a $\Bbb P^3$, hence
\[
  e_1+e_2+e_1'+e_2'+p=0.
\]
So
\[
  e+e'=-p.
\]
The converse is analogous.
\end{Proof}

{}From this lemma it follows that $f_1$ restricts to $D_1\iso E\times E$
as the quotient map to $E\times E/\iota'$, where $\iota'$ is the involution
$\iota'(p,e)=(p,-p-e)$. The curve $\Delta':=\{(-2e,e);\ e\in E\}$
is pointwise fixed under $\iota'$, while $\iota'$ acts as the
standard involution on the curve $(\Delta')^-:=\{(0,e);\ e\in E\}$.

Consider the change of coordinates (compare the proof of lemma \ref{L9}):
\[
  \psi_2\colon
  \left\{\begin{aligned}E\times  E &\to     E\times E\\
                        (p,e)      &\mapsto (p+e,-e).
  \end{aligned}\right.
\]

This maps $\Delta'$ to the diagonal $\Delta=\{(e,e);\ e\in  E\}$
and $(\Delta')^-$ to the antidiagonal $\Delta^-=\{(e,-e);\ e\in E\}$.
Moreover $\iota'$ becomes the involution $\tilde\iota$ interchanging the
two factors. From now on we shall identify $D_1$ with $E\times E$
via the isomorphism $\psi:=\psi_2\circ\psi_1$. Finally we denote by
$\bar\Delta$ the image of the diagonal $\Delta$ in
$S^2E=E\times E/\tilde\iota$.

\begin{Proposition}\label{P19}
The exceptional divisor $X\subset\Sigma_E^1$ intersects $D_1=E\times E$
in the diagonal $\Delta$.
\end{Proposition}

\begin{Proof}
The involution $\iota$ (resp.\ $\iota'$) is induced by a switching of cubic
scrolls in a singular quadric. A fixed scroll is precisely the unique
scroll containing $E$ in a rank 3 quadric.
\end{Proof}

\begin{Corollary}\label{C20}
$\cal F$ is singular along an elliptic scroll $S^2E$. The scroll $f_1(X)$
intersects $S^2E$ along $\bar\Delta$.
\end{Corollary}

\begin{Proposition}\label{P21}
The singular scroll $\pi_1(X)$ has degree 15. The curve $\Delta$ is mapped
4$:$1 to $E$ by $\pi_1$.
\end{Proposition}

\begin{Proof}
For the first part see \cite[prop.\ IV.4.7]{Hu}. The second statement
follows since $\Delta'\to E$ is given by $(-2e,e)\mapsto-2e$. It also follows
since the pencil of quadrics with vertex $p\in E$ contains 4 rank 3
quadrics (see the proof of proposition \ref{P16}).
\end{Proof}

We now turn our attention to the quintic hypersurface $V_2$.

\begin{Proposition}\label{P22}
Restricted to $D_1=E\times E$ the blowing down map $f_2$ is an isomorphism.
\end{Proposition}

\begin{Proof}
Fix some $e\in E$. Then the exceptional line in  $\Sigma_e^1$ and the curve
$E$ intersect transversally.
\end{Proof}

\begin{Proposition}\label{P23}
  \rom{(i)}
  The quintic hypersurface $V_2$ is ruled by an elliptic family of
  planes.

  \noindent\rom{(ii)}
  The map $\pi_2$ restricted to $D_1 \subset \widetilde V_2$ maps
  $D_1 \iso E {\times} E$ surjectively 2:1 onto a quintic elliptic
  scroll $S^2 E$. The scroll $S^2 E$ parametrizes those quadrics which
  are singular at a point of $E$. The quintic hypersurface $V_2$ is the
  trisecant scroll of $S^2 E$. It is singular exactly at $S^2 E$ (set
  theoretically).

  \noindent\rom{(iii)}
  Via $\pi_2$ the diagonal $\Delta \subset D_1 \iso E {\times} E$ is
  mapped to a degree 10 curve $\bar\Delta$ in $\PP^4$. The curve
  $\bar\Delta$ parametrizes the rank 3 quadrics through $E$.

  \noindent\rom{(iv)}
  The map $\pi_2$ gives an isomorphism of $\widetilde V_2 \sm D_1$ with
  $V_2 \sm S^2 E$.

  \noindent\rom{(v)}
  The rank of the differential of $\pi_2$ is 3 everywhere with the
  exception of $\Delta$ where it is 2.
\end{Proposition}

\begin{Proof}
  \rom{(i)}
  The fibre of $\widetilde V_2$ over a point $e \in E$ is mapped to the net of
  quadrics through the scroll $\Sigma_e^1$.

  \noindent\rom{(ii)}
  We have already seen that $f_1$ restricted to $D_1 \iso E {\times} E$
  factors through $S^2 E$. Hence using diagram \eqref{(D)} the same must
  be true for $\pi_2$. The map from $S^2E$ to $\PP^4$ given by $\pi_2$ is
  injective, which means that the image has degree at least 5. The ruling of
  $S^2E$ over a point $p\in E$ is mapped to the pencil of quadrics through $E$
  which are singular at $p$.  Now intersect $V_2$ with a general plane. Since
  $V_2$ is singular on the image of $S^2E$, this intersection is a plane curve
  with at least 5 singular points.  Since the map from
  $\widetilde V_2 \sm D_1$ to $V_2 \sm S^2 E$ is bijective, this curve
dominates
  the elliptic base curve of $\widetilde V_2$, and therefore, by the genus
formula,
  it cannot have more then 5 singular points. Thus $\pi_2(S^2E)$ has degree 5
  and, by the same argument, $V_2$ has no singularities outside $\pi_2(S^2E)$.
  If $C_0$ is a section of $S^2E$ with $C_0^2=1$ and $F$ is a fibre, then the
map
  from $S^2E$ to $\PP^4$ is given by the linear system $|C_0+2F|$. In fact,
  by $H_5-$invariance, the map is given by the complete linear system, in which
  case it is well known to be an embedding. We therefore identify $S^2E$ with
its
  image. It remains to show that $V_2$ is the trisecant scroll of $S^2E$.  Now,
  the curve $C_0$ moves in an elliptic family on $S^2E$, so each member is a
plane
  cubic curve.  Thus the planes of these curves are part of the trisecant
  scroll of $S^2E$.  Since $S^2E$ is the singular part of $V_2$, each such
  trisecant is contained in $V_2$ by Bezout. But the planes of the trisecant
  scroll cannot dominate the elliptic base curve, hence these planes must
coincide
  with the elliptic family of planes of $V_2$.

  \noindent\rom{(iii)}
  The curve $\Delta$ is the branch locus of the map
  $
    E {\times} E \to S^2 E \subset \PP^4
  $.
  It is well known that this is mapped to a curve $\bar\Delta$ of degree
  10 in $\PP^4$ (in fact the class of $\bar\Delta$ on $S^2 E$ is
  $4 C_0 - 2 F$ and the assertion follows from
  $
    (4 C_0 - 2 F) (C_0 + 2 F) = 10
  $).
  The assertion that $\bar\Delta$ parametrizes the rank 3 quadrics
  through $E$ follows from the description of the map $p_2$ and
  proposition \ref{P19}.

  \noindent\rom{(iv)}
  We have already seen that the map from $\widetilde V_2 \sm D_1$ to
  $V_2 \sm S^2 E$ is bijective. Since both sets are smooth, the claim
  follows.

  \noindent\rom{(v)}
  By \rom{(iv)} the rank of $d\pi_2$ is 3 outside $D_1$. Since the
  fibres of $\widetilde V_2$ are mapped to planes in $\PP^4$, it follows
  that the rank of the differential is at least 2 everywhere. Since
  $\Delta$ is the branch locus of the map $E {\times} E \to S^2 E$ the
  rank of $d \pi_2$ cannot be 3 along $\Delta$. It remains to prove
  that the rank of the differential is 3 on $D_1 \sm \Delta$. Let $x$ be
  a point on $D_1 \sm \Delta$ and let $E_x$ be the elliptic curve
  through $x$ which is mapped to a ruling of $S^2 E$. The differential
  of $\pi_2$ restricted to $E_x$ is 2 at $x$. Hence it is enough to see
  that the ruling $L_x = \pi_2(E_x)$ and the plane $\PP^2_x$ which is
  the image of the fibre of $\widetilde V_2$ containing $x$ meet
  transversally. For this it is enough to show that $L_x$ is not
  contained in $\PP^2_x$. But the intersection of $\PP_x^2$ with
  $S^2 E$ is a smooth plane cubic and does not contain a line.
\end{Proof}

\begin{Remark}\label{R24}
A general symmetric 5$\times$5 matrix with linear coefficients has
rank 3 along a curve of degree 20.
\end{Remark}
\bigskip

We are now in a position to connect the geometric approach of this paragraph
with the abstract approach from \ref{Par1}. To do this, recall the sections
$s_0,\ldots,s_4$ of $\cal L$ from proposition \ref{P6}.

\begin{Proposition}\label{P25}
There is an isomorphism $\widetilde V_2\iso\Bbb P_E^2$ such that the map
$\pi_2$ is given by $s_0,\ldots,s_4$.
\end{Proposition}

\begin{Proof}
The argument has two parts. First we identify $\widetilde V_2$ with ${\Bbb
P}(N_E(-2))$,
where $N_E(-2)$ is the twisted normal bundle of the elliptic curve $E\subset
{\PP^4}$.
Afterwards we show that ${\Bbb P}(N_E(-2))\iso {\Bbb P_E^2}$ and in fact also
the
existence of an $H_5$-isomorphism between ${\cal O}_{{\Bbb P}(N_E(-2))}(1)$ and
${\cal O}_{{\Bbb P}({\cal E})}(1)$.

We use the basis of $\mathrm H^\circ(I_E(2))$ given by
\[
  Q_i=x_i^2+ax_{i+2}x_{i+3}-\frac1ax_{i+1}x_{i+4}\quad(i\in\Bbb Z_5).
\]
The natural map
\[
  \mathrm H^\circ(I_E(2))\otimes\cal O_E\stackrel{\alpha}{\to}N_E^*(2)
\] is surjective and there is an exact sequence
\[
  0\to K\stackrel{\beta}{\to}\mathrm H^\circ(I_E(2))\otimes\cal O_E(-1)
  \stackrel{A}{\to}\mathrm H^\circ(I_E(2))\otimes\cal O_E
  \stackrel{\alpha}{\to}N_E^*(2)\to0
\]
with
\[
  A=\begin{pmatrix}     0 &  ax_4 &  -x_3 &   x_2 & -ax_1 \\
                    -ax_4 &     0 &  ax_2 &  -x_1 &   x_0 \\
                      x_3 & -ax_2 &     0 &  ax_0 &  -x_4 \\
                     -x_2 &   x_1 & -ax_0 &     0 &  ax_3 \\
                     ax_1 &  -x_0 &   x_4 & -ax_3 &     0
    \end{pmatrix}
\]
\medskip
(see \cite[p.\ 68]{Hu}). Dualising this sequence we get
\[
  \begin{CD} 
    \mathrm H^\circ(I_E(2))\spcheck\otimes\cal O_E
    @>(\tr A)=-A(1)>>
    \mathrm H^\circ(I_E(2))\spcheck\otimes\cal O_E(1)
    @>\tr\beta>>
    K^*
    @>>>
    0.
  \end{CD}
\]
Hence
\[
  K^*\iso N_E^*(3),
\]
i.e.,
\[
  K\iso N_E(-3).
\]

We want now to describe the map
\[
  \begin{CD}
    \Bbb P(N_E(-2))&\hookrightarrow&\Bbb P(\mathrm H^\circ(I_E(2)))\times E\\
                   &\searrow       &\downarrow \\
                   &               &\PP(\mathrm H^\circ(I_E(2)))
  \end{CD}
\]
where the horizontal map is given by the inclusion
\[
  N_E(-2)\stackrel{\beta(1)}{\hookrightarrow}\mathrm H^\circ(I_E(2))
  \otimes\cal O_E.
\]
We first want to identify the subbundle
\[
  \Bbb P(N_E(-2))=
  \{(p,Q);\ p\in E,\ Q\in\operatorname{Im}\beta(1)|_p\}
  \subset\Bbb P(\mathrm H^\circ(I_E(2))\times E.
\]

\begin{Claim}\label{P25-Claim1}
  $Q\in\operatorname{Im}\beta(1)|_p$ if and only if $Q$ contains the
  unique cubic scroll containing the secant $L_{\langle o,-p\rangle}$.
\end{Claim}

\begin{ProofwCaption}{Proof of the claim}
Consider the matrix
\[
  M'=\left(\frac{\partial Q_{3j}}{\partial x_i}\right)_{i,j}\quad
  (i,j\in\Bbb Z_5)
\]
from remark \ref{R17}. By proposition \ref{P16} \rom{(ii)} this has rank
$3$ on $E$. One easily checks that the entries of $A\tr M'$ are all
elements of $\mathrm H^\circ(I_E(2))$. Since $A$ has rank $2$ on $E$ the
sequence
\[
  \mathrm H^\circ(I_E(2))\otimes\cal O(-1)\stackrel{\tr M'}{\to}
  \mathrm H^\circ(I_E(2))\otimes\cal O_E\stackrel{A}{\to}
  \mathrm H^\circ(I_E(2))\otimes\cal O_E(1)
\]
is exact. Therefore
\[
  \operatorname{Im}\beta(1)=\operatorname{Im}\tr M'.
\]
Since there is a net of quadrics through a cubic scroll it suffices to
show that any quadric in the image of $\tr M'(p)$ contains the scroll.
For this it suffices to show that the secant lines $L_{\langle o,-p\rangle}$
and $L_{\langle\eta_5,\eta_{5-p}\rangle}$, where $\eta_5$ is a
non-zero 5-torsion point, are contained
in $Q$ (recall that the
cubic scroll in question is the union of all secants
$L_{\langle q,r\rangle}$ with $q+r=-p$). If $Q$ contains two secants in the
scroll it must contain the scroll by Bezout.

Now $\operatorname{Im}\tr M'(p)$ is spanned by the elements
\[
  M'_i(p)\begin{pmatrix} Q_0 \\
                         \vdots \\
                         Q_4
          \end{pmatrix}
\]
where $M'_i(p)$ is the $i$-th row of $M'$ evaluated at $p$. The
origin has coordinates $(0,a,-1,1,-a)$ and we can take $\eta_5$ to
be $(a,-1,1,-a,0)$. If $p$ has coordinates $(x_0,\ldots,x_4)$ then
$-p$ has coordinates $(x_0,x_4,x_3,x_2,x_1)$ and $-p-\eta_5$ has
coordinates $(x_4,x_3,x_2,x_1,x_0)$. Evaluating the quadrics in
$\operatorname {Im}\tr M'(p)$ on the secant lines one gets quadrics
in the coordinates $x_i$ which vanish on $E$. This proves the claim.
\noqed
\end{ProofwCaption}

 This shows that $\widetilde V_2=\Bbb P(N_E(-2))$ and that the map to
$\Bbb P^4$ is given by $\cal O_{\Bbb P(N_E(-2))}(1)$.

By \cite[proposition V.1.2]{Hu} the twisted normal bundle $N_E(-2)$ is
indecomposable with $c_1(N_E(-2))=\cal O_E(-1)$. By Atiyah's classification
\cite{At} $N_E(-2)\iso\cal E_E$ and in particular
$\Bbb P(N_E(-2))\iso \Bbb P^2_E$. Both bundles $N_E(-2)$ and $\cal E_E$
come with an $H_5$-action which covers the same action on $E$. Since $N_E(-2)$
resp.\ $\cal E_E$ are stable, and hence simple, the two $H_5$-actions on
$N_E(-2)$ and $\cal E_E$ differ at most by a character. But since the induced
actions on the respective determinants coincide, this character must be
trivial.
By construction $\cal L =\cal O_{\Bbb P(\cal E_E)}(1)$ (As a check note that
the representation of $H_5$ on both $\mathrm H^0(\cal L)$ (see proposition
\ref{P7}) and on
$
  \mathrm H^0(\cal O_{\Bbb P(\mathrm H^0(I_E(2)))} (1))
  =
  \mathrm H^0(I_E(2))\spcheck
$
are in each case derived from the Schr\"odinger representation by replacing
$\varepsilon$ by $\varepsilon^2$). In any case the above argument shows
that we have an $H_5$-isomorphism between $\cal O_{\Bbb P(N_E(-2))}(1)$ and
$\cal O_{\Bbb P(\cal E)}(1)$ and we are done.
\end{Proof}

\begin{Remark}\label{R27}
We have seen in proposition \ref{P23} (ii) that  $A_0=E\times E$ is
mapped 2$:$1 by $\pi_2$ onto an elliptic quintic scroll. Since the abelian
surfaces
$A_K$ and the bielliptic surfaces $B_{(a,b)}$ are numerically equivalent to
$A_0$ on $\widetilde V_2$ these surfaces must be mapped to surfaces in $\PP^4$
of degree 10.
\end{Remark}

\begin{Proposition}\label{P26}
There is an isomorphism $\widetilde V_1\iso\Bbb P^1_{S^2E}$ such that
$f_2$ is identified with $\rho$.
\end{Proposition}

\begin{Proof}
$f_2\colon\widetilde V_1\to\widetilde V_2$ is the blow up of $\widetilde V_2$
in the diagonal of $D_1$ after we have identified $D_1$ with $E\times E$
via the isomorphism $\psi$. Recall also that the map $\pi_2$ is
bijective outside $D_1$ and that $\pi_2$ restricted to $E\times E$ is a
2$:$1 branched covering onto its image whose branch locus is the
diagonal of $E\times E$. The map $\rho\colon\Bbb P^1_{S^2E}\to\Bbb P^2_E$ is
the blow up of $\Bbb P^2_E$ along the diagonal of $A_0=E\times E$. In view
of our identification of $\pi_2$ with the map given by $s_0,\ldots,s_4$ it
is enough to prove the following: The map $(s_0:\ldots:s_4)$ restricted to
$A_0=E\times E$ is a 2$:$1 branched covering with branch locus the
diagonal. But this is easy to see. Recall that the curve
$\{(q,q)\in E\times E\}\subset E\times\Bbb P^2$ is mapped 9$:$1 to the
antidiagonal $\{(q,-q);\ q \in E\}$
in $A_0$. By construction of $\cal L$ this shows that the degree of $\cal L$
restricted to the antidiagonal, and hence all its translates, is 2.
Moreover the degree of $\cal L$ restricted to $A_0$ is 10. It is
well known that then $A_0$ is mapped 2$:$1 onto a quintic elliptic scroll
with branch locus the diagonal (e.g. see \cite{HL}).
\end{Proof}

We are now ready to prove that the maps $\pi_1$ and $\pi_2$ give rise to
abelian and bielliptic surfaces of degree 15, resp.\ 10. Before we do
this, we recall from lemma 5 that
\[
  A_K \cap A_0 = \{(q,r) \in E \times E;\ 3 r + 2 q = 0\}
\]
resp.
\[
  B_{(a,b)} \cap A_0 = \{(q,r) \in E \times E;\ 3 r + 2 q = - \tau_{(a,b)}\}.
\]
We set
\[
  E_K = A_K \cap A_0, \quad E_{(a,b)} = B_{(a,b)} \cap A_0.
\]
Moreover, we consider the following curves on $E {\times} E$:
\[
  \Delta_p^- = \{(e,-e + p),\ e \in E\}.
\]
Under the quotient map $E {\times} E \to S^2 E$ these curves are mapped
to the rulings of the $\PP^1$-bundle $S^2 E$.

\begin{Lemma}\label{L27}
  The curves $E_K$, resp.\ $E_{(a,b)}$, intersect the curves
  $\Delta_p^-$ transversally in one point.
\end{Lemma}

\begin{Proof}
  A point $(e,-e + p)$ lies on $E_K$ if and only if $-e + 3 p = 0$,
  i.e., $e = 3 p$. Both curves are elliptic curves. Two curves on an
  abelian surface which do not coincide, meet transversally. The claim
  for the curve $E_{(a,b)}$ is proved in  exactly the same way.
\end{Proof}

\begin{Theorem}\label{T28}
  \rom{(i)}
  Let $A_K$ be a smooth element different from $A_0$ in the pencil
  $|-K|$ on $\PP^2_E \iso \widetilde V_2$. Then $\pi_2$ embeds $A_K$ as a
  smooth abelian surface of degree 10.

  \noindent\rom{(ii)}
  The bielliptic surfaces $B_{(a,b)}$ are also embedded as surfaces of
  degree 10 by the map $\pi_2$.
\end{Theorem}

\begin{Proof}
  \rom{(i)}
  By proposition \ref{P23} the map $\pi_2$ is an isomorphism outside
  $D_1 \iso A_0$. Hence it is sufficient to consider the intersection
  $E_K = A_K \cap A_0$. We first claim that $\pi_2$ restricted to  $E_K$
  is injective. The map $\pi_2$ identifies points $(q,r)$ and $(r,q)$.
  Assume that two such points lie on $E_K$. This implies that
  \[
    3 r + 2 q = 0, \quad 3 q + 2 r = 0.
  \]
  Subtracting these two equations from each other gives $q = r$, and
  hence $(q,r) = (r,q)$. Finally we have to check that the differential
  of $\pi_2$ restricted to $A_K$ is injective at the 25 points
  $E_K \cap \Delta$. For this recall that the kernel of $d\pi_2$ along
  $\Delta$ is given by the directions defined by the curves
  $\Delta_p^-$. The claim follows, therefore, from lemma \ref{L27}.
  The degree of the embedded surfaces is 10 by remark \ref{R27}.

  \noindent\rom{(ii)}
  The same proof goes through for the surfaces $B_{(a,b)}$.
\end{Proof}

\begin{Remarks}\label{R29}
  \rom{(i)}
  The pencil $|-K|$ contains 4 singular elements corresponding to the 4
  triangles in the Hesse pencil. These surfaces are mapped to
  translation scrolls of quintic elliptic curves where the translation
  parameter is a non-zero 3-torsion point.

  \noindent\rom{(ii)}
  All abelian surfaces $A_K$ are isogeneous to a product. Hence this
  construction does not give the general abelian surface in $\PP^4$. On
  the other hand we get all minimal bielliptic abelian surfaces in
  $\PP^4$ in this way (up to a change of coordinates).

  \noindent\rom{(iii)}
  The abelian surfaces in $\PP^4$ which are of the form
  $E {\times} F / \ZZ_3 {\times} \ZZ_3$ were studied by Barth and
  Moore in \cite{BM}. By their work the pencils which we have
  constructed above, are tangents to the rational sextic curve $C_6$ in
  the space of Horrocks-Mumford surfaces which parametrizes
  Horrocks-Mumford surfaces which are double structures on elliptic
  quintic scrolls.

  \noindent\rom{(iv)}
  The involution $\iota$ of proposition \ref{P8} induces the involution
  $x \mapsto -x$ on the surfaces $A_K$. The surfaces $B_{(a,b)}$ are
  identified pairwise, more precisely $\iota B_{(a,b)} = B_{(-a,-b)}$.
  This follows since the 8 characters $F_{(a,b)}$ are identified in this
  way by the Heisenberg involution on $\PP^2$.
\end{Remarks}

\begin{Theorem}\label{T30}
  \rom{(i)}
  Let $A_K$ be a smooth element different from $A_0$ in the pencil
  $|-K|$ on $\PP^2_E \iso \widetilde V_2$. Then the map $\pi_1$ embeds
  $\tilde A_K$ as a smooth non-minimal abelian surface of degree 15 in
  $\PP^4$.

  \noindent\rom{(ii)}
  The surfaces $\tilde B_{(a,b)}$ are embedded by $\pi_1$ as smooth
  bielliptic surfaces of degree 15.
\end{Theorem}

\begin{Proof}
  \rom{(i)}
  Again it is enough to look at the intersection of $\tilde A_K$ with
  $A_0$ on $\PP^1_{S^2 E} \iso \widetilde V_1$. The curves $\Delta_p^-$ are
  contracted by $\pi_1$. These are the only tangent directions which
  are in the kernel of the differential of $\pi_1$. Hence our claim
  follows again from lemma \ref{L27}. The double point formula reads
  \[
    d^2 = 10 d + 5 H K + K^2 - e.
  \]
  In our case $H K = 25$, $K^2 = -25$ and $e = 25$. This leads to
  \[
    d(d - 10) = 75
  \]
  and the only positive solution is $d = 15$.

  \noindent\rom{(ii)}
  The claim about the bielliptic surfaces can be proved in exactly the
  same way.
\end{Proof}

\begin{Remark}\label{R31}
  The degree of the surfaces $\tilde A_K$, resp.\ $\tilde B_{(a,b)}$,
  can also be computed by studying the linear system which maps
  $\widetilde V_1$ to $\PP^4$. We shall come back to this.
\end{Remark}
\bigskip

We shall now turn our attention to the quintic hypersurfaces which contain
the surfaces $A_K$, $B_{(a,b)}$, $\tilde A_K$ and $\tilde B_{(a,b)}$.

\begin{Proposition}\label{P33}
  The bielliptic surfaces $B_{(a,b)}$ lie on a unique quintic hypersurface,
  namely $V_2$.
\end{Proposition}

\begin{Proof}
  We consider the elliptic quintic scroll $S^2E\subset V_2\subset \PP^4$.
  Recall that $S^2 E$ is the quotient of $E \times E$ by the involution which
  interchanges the two factors. Let $C_{p_0}$ be the section of $S^2 E$
  which is the image of $\{p_0\} \times E$, resp.\ $E \times \{p_0\}$ in
  $S^2 E$ where $p_0$ is the origin of $E$ which we have chosen before.
  Note that the normal bundle of $C_{p_0}$ in $S^2 E$ is the degree 1 line
  bundle which is given by the origin.
  Let $F_{p_0}$ be the fibre over the origin of the map $S^2 E \to E$,
  $\{q_1,q_2\} \mapsto q_1 + q_2$. Moreover let $H$ be the hyperplane
  section of $S^2 E \subset \PP^4$. It follows immediately from our choice
  of the line bundle $\cal L$ in \ref{Par1} that
  \begin{equation}\label{(9)}
    H \sim C_{p_0} + 2 F_{p_0}.
  \end{equation}
  Next we consider the intersection $C_{(a,b)} = B_{(a,b)} \cap S^2 E$.
  The curve $C_{(a,b)}$ is by lemma \ref{L27} a section of $S^2E$. Now
  $E\times E\to S^2E$ is ramified along the diagonal and maps the curve
  $E_{(a,b)}$ isomorphically to $C_{(a,b)}$ so, combining with lemma
  \ref{L5}\rom{(ii)}, the intersection of $C_{(a,b)}$ with the diagonal is
  twice the set $\{(p,p);\ 5p=-\tau_{(a,b)}\}$. Thus $C_{(a,b)} \equiv
C_{p_0}+12F$.
  In fact the intersection with the diagonal goes by the map $S^2E\to E$ to
  $\{4p;\ 5p=-\tau_{(a,b)}\}$, which summed up is $\tau_{(a,b)}$.  Therefore
  \begin{equation}\label{(10)}
    C_{(a,b)} \sim C_{p_0} + 11 F_{p_0} + F_{(a,b)}
  \end{equation}
  where  $F_{(a,b)}$ is the fibre over the 3-torsion point $\tau_{(a,b)}$.
  Finally recall that
  \begin{equation}\label{(11)}
    K \sim - 2 C_{p_0} + F_{p_0}.
  \end{equation}
  Let $Q$ be a quintic containing $B_{(a,b)}$. We first claim that $Q$ must
  contain $S^2 E$. In order to  see this look at the exact sequence
  \[
    0
    \to
    \cal O_{S^2 E}(5 H - C_{(a,b)})
    \to
    \cal O_{S^2 E}(5 H)
    \to
    \cal O_{C_{(a,b)}}(5 H)
    \to
    0.
  \]
  It follows from formulas \eqref{(9)}, \eqref{(10)}, and \eqref{(11)} that
  \begin{equation}\label{(12)}
    5 H - C_{(a,b)}
    \sim
    4 C_{p_0} - F_{p_0} - F_{(a,b)}
    \sim
    -2 K + (F_{p_0} - F_{(a,b)}).
  \end{equation}
  It is well known that $h^0(\cal O_{S^2 E}(-2 K + F_p - F_q)) = 0$
  unless $2p = 2q$ (cf. \cite{CC}). Hence $Q$ must contain $S^2 E$. We next
claim
  that $Q = V_2$. In order to see this, consider a plane on $V_2$, i.e., a
  trisecant plane of $S^2 E$. Both $S^2 E$ and $B_{(a,b)}$ intersect such a
  plane in different irreducible cubic curves. Since $Q$ has degree 5 it
  must contain this plane and hence $V_2$. By reasons of degree this
  implies $Q = V_2$.
\end{Proof}

\begin{Remark}\label{R34}
  The surfaces $A_K$ lie on 3 independent quintics. This follows
  e.g.\ from the fact that $A_K$ is the zero-scheme of a section $s$ of
  the Horrocks-Mumford bundle $F$ (where we normalize $F$ such that $c_1(F) =
5$,
  $c_2(F) = 10$). In other words there is an exact sequence
  \[
    0
    \to
    \cal O_{\PP_4}
    \stackrel{s}{\to}
    F
    \to
    I_{A_K}(5)
    \to
    0.
  \]
  The claim follows from $h^\circ(F) = 4$ (see \cite{HM}).
\end{Remark}

\medskip

Note that if we replace $B_{(a,b)}$ by a surface $A_K$ in the proof of
proposition \ref{P33} we obtain $-2 K$ in formula \eqref{(11)}. Since
$
  h^\circ(\cal O_{S^2 E}(-2 K)) = 2
$,
this gives rise to two more elements in
$
  \mathrm H^\circ(\cal O_{S^2 E}(-2 K))
$
which vanish along $C_K = A_K \cap S^2 E$.
These can be lifted to $H_5$-invariant quintics in $\PP^4$. It is then easy
to check that any $H_5$-invariant quintic which contains $C_K$ must contain
$A_K$. (Look at the intersection of the quintic with the cubic curves on
$A_K$. Unless the quintic contains these curves, this would split up into
two $H_5$-orbits of length 9 and 6 resp., a contradiction).
In this way one can also prove the existence of 3 independent quintics
through the surfaces $A_K$.

We now turn our attention to the degree 15 surfaces.

\begin{Proposition}\label{P35}
  \rom{(i)} The non-minimal abelian surfaces $\tilde A_K$ lie on exactly
  three quintic hypersurfaces. They are linked $(5,5)$ to translation
  scrolls.
  \par\noindent\rom{(ii)} The non-minimal bielliptic surfaces
  $\tilde B_{(a,b)}$ lie on a unique quintic hypersurface, namely $V_1$.
\end{Proposition}

\begin{Proof}
  \rom{(i)} Let $\bar H_1$ be the hyperplane section on $\PP^1_{S^2 E}$
  given by the map $\pi_1\colon \PP^1_{S^2 E} \to V_1 \subset \PP^4$. Let
  $\bar C_{p_0}$, resp.\ $\bar F_{p_0}$ be the fibres over the curves
  $C_{p_0}$, resp.\ $F_{p_0}$ in $S^2 E$ with respect to the map
  $\PP^1_{S^2 E} \to S^2 E$. The classes $\bar H_1$, $\bar C_{p_0}$ and
  $\bar F_{p_0}$ generate the Neron-Severi group of $\PP^1_{S^2 E}$.
  Under $\pi_1$ the surface $\bar F_{p_0}$ is mapped isomorphically to a
  cubic scroll, while $\bar C_{p_0}$ is mapped birationally to the cone
  over an elliptic curve of degree 4 in $\PP^3$ with vertex in the
  origin $p_0$ of the elliptic curve in $\PP^4$ for which $V_1$ is the
  secant variety. Thus we get the following intersection numbers:
  $\bar H_1^3 = 5$,
  $\bar H_1^2 \bar C_{p_0} = 4$,
  $\bar H_1^2 \bar F_{p_0} = 3$,
  $\bar H_1 \bar C_{p_0} \bar F_{p_0} = \bar H_1 \bar C_{p_0}^2 = 1$
  and
  $
    \bar C_{p_0}^2 \bar F_{p_0}
    =
    \bar C_{p_0} \bar F_{p_0}^2
    =
    \bar H_1 \bar F_{p_0}^2
    =
    \bar F_{p_0}^3
    =
    \bar C_{p_0}^3
    =
    0
  $.
  The exceptional divisor $X = E_\Delta$ is a section of the $\PP^1$-bundle
  $\PP^1_{S^2 E}$. Under $\pi_1$ it is mapped birationally onto a ruled
  surface of degree 15. From this information it is easy to compute that
  numerically $X \equiv \bar H_1 - 2 \bar C_{p_0} + 6 \bar F_{p_0}$. It is
  also straightforward to check that for the canonical divisor
  $K \equiv -2 \bar H_1 + \bar C_{p_0} + 2 \bar F_{p_0}$. In fact we claim
  that these equalities are also true with respect to linear equivalence.
  It is enough to check this on a section of the composite
  projection $\PP^1_{S^2 E} \to S^2 E \to E$.
  We consider the curve $D = X \cap D_1 = X \cap (E \times E)$. On
  $E \times E$ this is the diagonal by proposition \ref{P19}. Using the
  projection $\PP^1_{S^2 E} \to S^2 E$ we can identify the section $X$ with
  $S^2 E$. Then by \cite[lemma IV.4.4]{Hu} we have
  $D \sim C_{p_0} + 12 f_{p_0}$ on $X$. Since $X$ is mapped to a ruled
  surface of degree 15 we have $\bar H_1|_X \equiv C_{p_0} + 7 F_{p_0}$.
  Since $D$ is mapped by multiplication with $-2$ four to one onto the
  elliptic curve $E$ and since $E$ is embedded by $|5 p_0|$ we have in fact
  that this equality also holds with respect to linear equivalence. We also
  note that $X$ restricted to $D$ is trivial. This follows since $X$
  restricted to $A_0$ is $D$ and since the normal bundle of $D$ in $A_0$ is
  trivial. Hence in order to check that
  $
    X \sim \bar H_1 - 2 \bar C_{p_0} + 6 \bar F_{p_0}
  $
  it is enough to prove that the restriction of
  $\bar H_1 - 2 \bar C_{p_0} + 6 \bar F_{p_0}$ to $D$ is trivial. This
  follows from
  \[
    (-C_{p_0} + 13 F_{p_0}) (C_{p_0} + 12 p_0) \sim 0
  \]
  which has to be read as an equality of divisors on $D \iso E$.
  Next we want to prove that
  \[
    K \sim -2 \bar H_1 + \bar C_{p_0} + 2 \bar F_{p_0}.
  \]
  We know that $(K + X)|_D \sim K_X|_D \sim -25 p_0$. The first is the
  adjunction formula, the second follows from
  \[
    (-2 C_{p_0} + F_{p_0}) (C_{p_0} + 12 F_{p_0}) \sim -25 p_0.
  \]
  Since $X|_D \sim 0$ it is now enough to show that
  $-2 \bar H_1 + \bar C_{p_0} + 2 \bar F_{p_0}$ restricted to $D$ is
  linearly equivalent to $-25 p_0$. This follows from
  \[
    (-C_{p_0} - 12 F_{p_0}) ( C_{p_0} + 12 F_{p_0}) \sim -25 p_0.
  \]
  Hence we have proved that
  \begin{equation}\label{(12a)}
    \bar X \sim \bar H_1 - 2 \bar C_{p_0} + 6 \bar F_{p_0},
    \quad
    K \sim -2 \bar H_1 + \bar C_{p_0} + 2 \bar F_{p_0}.
  \end{equation}
  Since $A_K$ is anticanonical on $\PP^2_E$ we get that
  $
    \tilde A_K
    \sim
    -K + X
    \sim 3
    \bar H_1 - 3 \bar C_{p_0} + 4 \bar F_{p_0}
  $.
  Furthermore $-K \sim D_1$, and $D_1$ is contracted under $\pi_1$ to the
  curve $E$. Since twice the anticanonical divisor on  $S^2 E$ moves in a
  pencil it follows that also $4 \bar C_{p_0} - 2 \bar F_{p_0}$ moves in a
  pencil on $\PP^1_{S^2 E}$. Note that
  \begin{equation}\label{(13)}
    5 \bar H_1
    \sim
    -K + \tilde A_K + (4 \bar C_{p_0} - 2 \bar F_{p_0}).
  \end{equation}
  The members of the pencil $(4 \bar C_{p_0} - 2 \bar F_{p_0})$ on
  $\PP^1_{S^2 E}$ are mapped to the translation scrolls of $E$. Take such a
  translation scroll (which is not a quintic elliptic scroll). Then it is a
  Horrocks-Mumford surface (cf. \cite{Hu2}, \cite{BHM}) and hence lies on three
  quintics of which  $V_1$ is one. We can choose a pencil of quintics through
  such a scroll which does not contain $V_1$. All these quintics contain $E$.
  They cut out a pencil of residual surfaces and it follows from \eqref{(13)}
  that this is just the pencil formed by the surfaces $\tilde A_K$. Hence
  every such surface is linked $(5,5)$ to a translation scroll $S$. Now
  consider the well known liaison sequence (cf.\ \cite{PS}):
  \[
    0
    \to
    I_{S \cup \tilde A_K}(5)
    \to
    I_{\tilde A_K}(5)
    \to
    \omega_S
    \to
    0.
  \]
  We have $h^\circ(I_{S \cup \tilde A_K}(5)) = 2$,
  $h^1(I_{S \cup \tilde A_K}(5)) = 0$. Since $S$ is a a Horrocks-Mumford
  surface $\omega_S = \cal O_S$ and hence $h^\circ(\omega_S) = 1$. It follows
that
  $h^\circ(I_{\tilde A_K}(5)) = 3$.

  \noindent\rom{(ii)}
  Now consider a bielliptic surface $\tilde B_{(a,b)}$. If
  $\tilde B_{(a,b)}$ lies on two quintics, it would be linked to a surface
  $T$ in the numerical equivalence class of
  $4 \bar C_{p_0} - 2 \bar F_{p_0}$. Since $3 \tilde B_{(a,b)}$ is linearly
  equivalent to $3 \tilde A_K$, we must have that $3 T$ is linearly
  equivalent to $3 (4 \bar C_{p_0} - 2 \bar F_{p_0})$ while $T$ is not
  linearly equivalent to $4 \bar C_{p_0} - 2 \bar F_{p_0}$. But in this
  numerical equivalence class the only effective divisors are the pencil
  $(4 \bar C_{p_0} - 2 \bar F_{p_0})$ and the three divisors
  $4 \bar C_{p} - \bar F_{p_0} - \bar F_\tau$ where $\tau$ is a non-trivial
  2-torsion point.
\end{Proof}

\begin{Remark}\label{R36}
  At this point we would like to say a few more words about liaison. As said
  before, the space $\Gamma(F)$ of sections of the Horrocks-Mumford bundle has
  dimension 4. The three-dimensional space $\PP \Gamma = \PP(\Gamma(F))$
  parametrizes the Horrocks-Mumford (\rom{HM}) surfaces $X_s = \{s = 0\}$
  where $0 \ne s \in \Gamma(F)$. The space of Heisenberg
  invariant quintics is related to $\Gamma(F)$ via the isomorphism
  \[
    \Lambda^2 \Gamma(F) \iso \Gamma_{H_5}(\cal O_{\PP^4}(5))
  \]
  given by the natural map
  $
    \Lambda^2 \Gamma(F)
    \to
    \Gamma(\Lambda^2 F)
    =
    \Gamma(\cal O_{\PP^4}(5))
  $.
  Set
  \[
    \PP_{H_5}^5
    =
    \PP(\Gamma_{H_5}(\cal O_{\PP^4}(5)))
    \iso
    \PP(\Lambda^2 \Gamma(F)).
  \]
  In $\PP_{H_5}^5$ we consider the Pl\"ucker quadric $G = G(1,3)$ of
  decomposable tensors.

  If $X_{s_0} = \{s_0 = 0\}$ is a \rom{HM}-surface, then
  \[
    \Gamma(I_{X_{s_0}}(5)) = \{s_0 \wedge s;\ s \in \Gamma(F)\}.
  \]
  This defines a $\PP^2$ of decomposable tensors. In $G = G(1,3)$ this is an
  $\alpha$-plane, i.e., a plane of lines through one point. In this way we
  get a bijection between $\PP\Gamma$ and the set of all $\alpha$-planes in
  $G$.
  Now consider a line in an $\alpha$-plane, i.e., a pencil of quintics
  spanned by quintics of the form $s_0 \wedge s_1$ and $s_0 \wedge s_2$. This
  give rise to a complete intersection
  \[
    Y_1 \cap Y_2 = X_{s_0} \cup X'
  \]
  where $X'$ is of degree 15. By the liaison sequence which we have already
  used before, we find that $h^\circ(I_{X'}(5)) = 3$. The space of quintics
  is spanned by $s_0 \wedge s_1$, $s_0 \wedge s_2$ and $s_1 \wedge s_2$. To
  see that $s_1 \wedge s_2$ is contained in this space consider
  $X' \sm X_{s_0}$. At these points $s_0$ does not vanish and $s_1$ and $s_2$
  are linearly dependent of $s_0$. It follows that $s_1 \wedge s_2 = 0$.
  Hence all quintics through $X'$ are in particular $H_5$-invariant and the
  space of these quintics is a $\beta$-plane in $G$, i.e., a $\PP^2$ of lines
  which lie in a fixed plane in $\PP^3$.

  Let $\tilde A_K$ be one of our non-minimal degree 15 abelian surfaces.
  By proposition \ref{P35} the surface $\tilde A_K$ is linked to a
  translation scroll. Hence it defines a $\beta$-plane in $G$. For every line
  in this plane there exists exactly one $\alpha$-plane intersecting the
  $\beta$-plane in this line. Hence every such line gives rise to liaison
  with an \rom{HM}-surface. In this way $\tilde A_K$ is linked to a
  2-dimensional family of \rom{HM}-surfaces. This 2-dimensional family is
  parametrized by a linear $\PP^2$ in $\PP \Gamma$. Since the singular
  \rom{HM}-surfaces form an irreducible surface of degree 10 \cite{BM} in
  $\PP \Gamma$, it follows in particular that $\tilde A_K$ is linked to
  smooth abelian surfaces.
\end{Remark}

\medskip

We want to conclude this paragraph with a short discussion of the
6-secants of the surfaces $A_K$ and $B_{(a,b)}$. The 6-secants of the
surfaces $A_K$ are exactly the 25 Horrocks-Mumford lines \cite{HM}.

The 6-secant formula from \cite{L} shows that the surfaces $B_{(a,b)}$
either also have 25 6-secants or infinitely many. In fact there are exactly
25 6-secants and we shall now describe them. First note that every 6-secant
of $B_{(a,b)}$ must lie in one of the planes of $V_2$ (it must be contained
in $V_2$ by reasons of degree and since the base of the bundle $\PP^2_E$ is
elliptic,  it must be in one of the fibres). Now fix a point $e \in E$ and
let $f=f_{2}\! \! \shortmid_{\Sigma^{1}_{e}}$ be the blowing down map
$\Sigma_e^1 \to \PP_e^2$. If we interpret
$\Sigma_e^1$ as a cubic scroll in $\PP^4$, then $\Sigma_e^1$ consists of
the secants of the elliptic quintic curve $E$ joining points $p$ and $q$ with
$p + q = e$. We denote this curve by $E_e \subset \Sigma_e^1$. The map
$f\colon \Sigma_e^1 \to \PP_e^2 \subset \PP^4$ is given by the linear
system $|X_e + l|$ where $l$ are the fibres of the $\PP^1$-bundle
$\Sigma_e^1$ and  $X_e$ is the exceptional line. By what we have said
before, we have $l \cap E_e = \{p,q\}$ with $p + q = e$. One also computes
easily that $X_e \cap E_e = \{-2 e\}$ where we consider $-2 e$ as a point
on $E_e$. It follows that the map $f$ is given by the linear system
$(-2 e) + e + p_0 \sim 5 p_0 - 2 e_0$ on $E_e$ where $2 e_0 = e$.
The composite $\pi_{2} \circ f_{2}$ maps both curves $E_{e}$ and
$\iota' E_{e} \subset E\times E$ onto the same plane cubic curve in $\PP^{4}$.
When considering the change of coordinates $\psi_2$ followed by the involution
interchanging the factors on $E\times E$, $E_{e}$ is mapped to the first
factor,
while $\iota 'E_{e}$ is mapped to the second factor.  Therefore, by lemma
\ref{L5} (ii), the intersection $\tilde B_{(a,b)} \cap \iota' (E_{e})$ consists
of the points $\{(p,-p-e);\ 2p=e-\tau_{(a,b)}\}$. These have the same image on
the scroll as the points $p$ on $E_{e}$ with $2p=e-\tau_{(a,b)}$.
These are the points  $p_i = e_0 + \tau_{(a,b)} + \tau_i$, where
the $\tau_i$, $i=0,1,2,3$, are the 2-torsion points with $\tau_0 = p_0$.
The points $f(p_i)$, $i=1,2,3$, are collinear if and only if
$\sum_{i = 1}^3 p_i + 2 e_0 = p_0$ which is equivalent to $5 e_0 = p_0$,
i.e., $e_0$ is a 5-torsion point on $E$. In this case we also get
$e_0 = -4 e_0 = -2 e \in X_e \cap E_e$. Hence for $e$ a 5-torsion point on
$E$ the line through these 3 points intersects $B_{(a,b)}$ in addition in
3 points on its plane cubic in $\PP_e^2$, hence is a 6-secant. The above
discussion also shows that there are only finitely many 6-secants, and
hence we have found them all.

\section{Cremona transformations}\label{Par3}

In this paragraph we want to explain how the abelian, resp.\ bielliptic
surfaces of degree 15 can be constructed via Cremona transformations.

It is known that the quadrics through an elliptic quintic curve
$E \subset \PP^4$ define a Cremona transformation
$\Phi\colon \PP^4 \dashto \PP^4$  \cite{Sem}. Via $\Phi$ the secant variety
of $E$ is mapped to a quintic elliptic scroll $S^2 E$ in $\PP^4$. The
exceptional
locus of $\Phi$ is mapped to the trisecant variety of this scroll. The
cubics through the quintic elliptic scroll $S^2 E$ define the inverse of
$\Phi$. Under $\Phi^{-1}$ the trisecant variety of $S^2 E$ is mapped to
$E$, while the exceptional divisor is mapped to the secant variety of $E$.
We refer to \cite{Sem} for details on the geometry of this transformation.

We consider the scroll $S^2 E$ whose trisecant variety equals $V_2$. On
$S^2 E$ there exist three elliptic 2-sections $E_i$, $i=1,2,3$, such that
$S^2 E$ is the translation scroll of $E_i$ defined by a non-zero 2-torsion
point $p_i$. Let $\widetilde V_2 \iso \PP^2_E$ be the desingularization of
$V_2$. Then the curves
$
  \Delta_i = \{(e + p_i,e);\ e \in E\} \subset A_0 \iso E \times E
$
are mapped 2$:$1 to the elliptic curves $E_i$. Let
$\Delta_E \subset A_0 \iso E \times E$ be the diagonal. Let $\bar H_2$ be the
line bundle on $\widetilde V_2$ given by the map
$\pi_2\colon \widetilde V_2 \to V_2 \subset \PP^4$.

Next consider the isomorphism
\[
  \begin{aligned}
    \phi_i\colon E \times \PP^2 &\to     E \times \PP^2 \\
                 (e,x)          &\mapsto (e + p_i,x).
  \end{aligned}
\]
This commutes with the diagonal action of $\ZZ_3 \times \ZZ_3$ on
$E \times \PP^2$ and hence induces an isomorphism
\[
  \tilde\phi_i\colon \widetilde V_2 \to \widetilde V_2
\]
which maps $\Delta_E$ to $\Delta_i$. Note also that $\tilde\phi_i$ maps the
surfaces $A_K$, resp.\ $B_{(a,b)}$ to themselves. Recall that $\widetilde V_1$
is the blow-up of $\widetilde V_2$ along $\Delta_E$. Let $\widetilde V_1^{(i)}$
be
the blow-up of $\widetilde V_2$ along $\Delta_i$. Then $\tilde\phi_i$ induces
an isomorphism
\[
  \bar\phi_i\colon \widetilde V_1 \to \widetilde V_1^{(i)}
\]
such that the diagram
\[
  \begin{CD}
    \widetilde V_1 @>\bar\phi_i>>   \widetilde V_1^{(i)} \\
    @V\rho VV                   @VV\rho^{(i)}V   \\
    \widetilde V_2 @>\tilde\phi_i>> \widetilde V_2
  \end{CD}
\]
commutes where the vertical maps are the blowing down maps. $X$ is the
exceptional locus of $\rho$. Let $X_i$ be the exceptional loci of the maps
$\rho^{(i)}$. By abuse of notation we denote the pullback of
$\bar H_2$ by $\rho$ also by $\bar H_2$. We denote the pullback of
$\bar H_2$ by $\rho^{(i)}$ by $\bar H_2^{(i)}$. The Cremona transformation
defined by the quadrics through $E_i$ gives rise to the linear system
$|2\bar H_2^{(i)} - X_i|$ on $\widetilde V_1^{(i)}$. Note that
\[
  (\bar\phi_i)^* (2\bar H_2^{(i)} - X_i)
  \sim
  \rho^*(\tilde\phi_i)^*(2 \bar H_2) - X
  \sim
  2 \bar H_2 - X.
\]
The latter follows since $\bar H_2$ restricted to $\Delta_E$ has degree 10
and this implies that translation by a 2-torsion point leaves the linear
equivalence class invariant.

\begin{Proposition}\label{P37}
  $\bar H_1 \sim 2 \bar H_2 - X$.
\end{Proposition}

\begin{Proof}
  We first claim that $\bar H_1 \equiv 2 \bar H_2 - X$. The N\'eron-Severi
  group of $\widetilde V_1$ is generated by $\bar H_2$, $X$ and $\Sigma^1$
  where $\Sigma^1$ denotes the class of a fibre of $\widetilde V_1 \to E$.
  Restriction to such a fibre implies immediately that
  $\bar H_1 \equiv \alpha \Sigma^1 + 2 \bar H_2 - X$.
  To compute $\alpha$ we use $\bar H_1^3 = 5$. Since $(\Sigma^1)^2 = 0$
  this implies
  \begin{equation}\label{(14)}
    5 = 9 \alpha + (2 \bar H_2 - X)^3.
  \end{equation}

  Now $\bar H_2^2\cdot X=0$, and $\bar H_2\cdot X^2=-10$, since $X$
  is blown down to the diagonal $\Delta_E\subset A_0$ and $\Delta_E\cdot\bar
H_2=10$.
  On the other hand $X^3=-25$ from our computations in the proof of
  proposition \ref{P35}, thus $\alpha =0$ as claimed. In order to prove the
  proposition it is now enough to consider the restriction of $\bar H_1$,
  resp.\ $2 \bar H_2 - X$ to the section $D$ which we have already used in
  the proof of proposition \ref{P35}. Via $\rho$ the curves $D$ and
  $\Delta_E$ are identified. We know that $\bar H_1$ restricted to $D$ is
  linearly equivalent to $20 p_0$. On the other hand $\bar H_2$ restricted
  to $\Delta_E$ is linearly equivalent to $10 p_0$ (this can be seen e.g.\ by
  using proposition \ref{P19} and the construction of the line bundle
  $\cal L$). We also have already seen in the proof of proposition
  \ref{P35} that the restriction of $X$ to $D$ is trivial. This proves the
  proposition.
\end{Proof}

\begin{Corollary}\label{C38}
  The map $\pi_1\colon \widetilde V_1 \to V_1$ is given by the complete linear
  system $|\bar H_1| = |2 \bar H_2 - X|$.
\end{Corollary}

\begin{Proof}
  We have to show that the (affine) dimension of the linear system
  $|2 \bar H_2 - X|$ is five. We consider the exact sequence
  \[
    0
    \to
    \cal O_{\widetilde V_1}(2 \bar H_2 - A_0 - X)
    \to
    \cal O_{\widetilde V_1}(2 \bar H_2 - X)
    \to
    \cal O_{A_0}(2 \bar H_2 - \Delta_E)
    \to
    0.
  \]
   Since $(2H_2-A_0-X)\cdot H_2\cdot\Sigma^1 =-1$ it follows that
 $h^\circ( \cal O_{\widetilde V_1}(2 \bar H_2 - A_0 - X)) =0$. Hence we have
 the inclusion
   \[
    0
    \to
    \mathrm H^\circ(\cal O_{\widetilde V_1}(2 \bar H_2 - X))
    \to
    \mathrm H^\circ(\cal O_{A_0}(2 \bar H_2 - \Delta_E)).
  \]
  The linear system $|2 \bar H_2 - \Delta_E|$ restricted to $A_0$ has degree
  0 on the curves $\{(p,-p + e);\ p \in E\}$ and degree 5 on the curves
  $\{(e,p);\ p \in E\}$. It follows that
  $h^\circ(\cal O_{A_0}(2 \bar H_2 - \Delta_E)) = 5$ and this proves the
  corollary.
\end{Proof}

\begin{Corollary}\label{C39}
  The hyperplane bundle of $\tilde A_K$ is of the form
  $2 H' - \sum_{i = 1}^{25} E_i$ where $H'$ is a polarization of type
  $(1,5)$ on the minimal model of $\tilde A_K$.
\end{Corollary}

\begin{Proof}
  Immediately from proposition \ref{P37}.
\end{Proof}

\end{document}